\documentclass{article}
\usepackage{amsmath,amssymb,amsfonts}
\usepackage{graphicx}
\usepackage{color}

\usepackage[all]{xy}

\def\rank{\mathrm{rank}}
\def\deg {\mathrm{deg}}

\def\Jac{\mathrm{Jac}}

\def\inv{\mathrm{inv}}

\newtheorem{theorem}{Theorem}
\newtheorem{proposition}{Proposition}
\newtheorem{corollary}{Corollary}
\newtheorem{lemma}{Lemma}

\newtheorem{example}{Example}
\newtheorem{definition}{Definition}
\newenvironment{proof}[1][Proof]{\noindent\textit{#1.} }{\hfill$\Box$\medskip}

\title{Two-valued groups, Kummer varieties and integrable billiards}
\author{V. M. Buchstaber (*), V. Dragovi\' c (**)}

\date{}

\begin{document}

\maketitle

\medskip
\centerline{(*) Steklov Mathematical Institute, Russian Academy of
Sciences} \centerline{ Gubkina Street 8, 119991 Moscow, Russia}

\smallskip

\centerline{e-mail: {\tt buchstab@mi.ras.ru}}

\medskip
\centerline{(**)Mathematical Institute SANU}

\centerline{Kneza Mihaila 36, 11000 Belgrade, Serbia}

\smallskip

\centerline{Mathematical Physics Group, University of Lisbon}

\smallskip

\centerline{e-mail: {\tt vladad@mi.sanu.ac.rs}}

\

\begin{abstract}
\smallskip
A natural and important question of study two-valued groups
associated with hyperelliptic Jacobians and their relationship with
integrable systems is motivated by seminal examples of relationship
between algebraic two-valued groups related to elliptic curves and
integrable systems such as elliptic billiards and celebrated
Kowalevski top. The present paper is devoted to the case of genus 2,
to the investigation of algebraic two-valued group structures on
Kummer varieties.

One of our approaches is based on the theory of $\sigma$-functions.
It enables us to study the dependence of parameters of the curves,
including rational limits. Following this line, we are introducing a
notion of $n$-groupoid as natural multivalued analogue of the notion
of topological groupoid.

Our second  approach is geometric. It is based on a geometric
approach to addition laws on hyperelliptic Jacobians and on a recent
notion of billiard algebra. Especially important is connection with
 integrable billiard systems within confocal quadrics.

The third approach is based on the realization of the Kummer variety
in the framework of moduli of semi-stable bundles, after Narasimhan
and Ramanan. This  construction of the two-valued structure is
remarkably similar to the historically first example of topological
formal two-valued group from 1971, with a significant difference:
the resulting bundles in the 1971 case  were "virtual", while in the
present case the resulting bundles are effectively realizable.

\end{abstract}



\newpage

\tableofcontents

\newpage

\section{Introduction}\label{sec:intro}

\medskip

The study of structure   of multivalued groups has been started in
1971 (see \cite{BN}) within the study of characteristic classes of
vector bundles, at that time as the structure of formal (local)
multivalued groups. In 1990, in \cite{Buc0} has been introduced
algebraic two-valued group structure on $\mathbb {C}P^1$ by using
addition theorems on elliptic curves. The structure of algebraic
multivalued groups has been studied since then extensively (see
\cite{Buc} and references therein) in different contexts.

The deep connection between two-valued groups and integrable systems
is well-known, see \cite {BV}. However, it has got a new impetus
with a recent paper \cite {Drag2}, where a connection between
two-valued groups on $\mathbb CP^1$  and the celebrated Kowalevski
top  has been discovered. There, it was shown that the "mysterious"
Kowalevski change of variables (see, for example, \cite {Aud} for
the terminology) corresponds to the operation in two-valued group
$W/\tau$, where $W$ is an elliptic curve, in the standard
Weierstrass model of elliptic curve, with $\tau$ as the canonical
involution. More detailed description of this two-valued group is
presented in Section \ref{subsection:defexample}, see Example
\ref{ex:eplus}. This deep connection of Kowalevski integration
procedure with a structure of elliptic curve, on the first glance,
may be in surprising contrast to a well known fact that the
integration of the Kowalevski top is performed on Jacobian of a
genus two curve. Moreover, as it has been shown in \cite {Drag2},
the Kowalevski integration procedure can be interpreted as a certain
deformation of the two-valued group $p_2$ (see Example
\ref{ex:cplus}) toward the structure on $W/\tau$. The structure of
$p_2$  is the  rational limit of the one on $W/\tau$. In addition,
in \cite {Drag2} it has been proven the equivalence between
associativity of the two-valued group on $W/\tau$ with the Poncelet
theorem of pencils of conics, in a case of triangles.

Higher genera Poncelet type problems, integrable billiards and
so-called billiard algebra associated with pencils of quadrics in
arbitrary dimensions have been studied in \cite {DR2} and \cite
{DR}, where the relationship with hyperelliptic curves and their
Jacobians has been considered (see also \cite {DR3}).

From our previous experience and work with elliptic case, as a
natural and important question arises the study two-valued groups
associated with hyperelliptic Jacobians and their relationship with
integrable systems. Especially important is connection with
integrable billiards within pencils of quadrics.

The present paper is devoted to the case of genus 2. In this case
one comes to the problem of investigation of two-valued group
structures on Kummer varieties.

We describe the structure of algebraic two-valued groups on Kummer
varieties, its rational limits and relationship with integrable
systems. The Kummer variety is a classical and well-studied object
of the algebraic geometry. It appears as a variety of orbits of a
Jacobian of a curve of genus 2, factorized by the hyperelliptic
involutive automorphism. Thus, there are various algebro-geometric
and analytical structures on it, related to the moduli of curves of
genus 2. We want to understand and present the structure of an
algebraic two-valued group associated to different appearances of
the Kummer varieties. Let us mention that from the point of view of
differential geometry and the theory of abelian functions of genus
2, Kummer varieties have been studied by Baker.

One side of our approach is based on the theory of
$\sigma$-functions. The choice of $\sigma$ functions is motivated by
our wish to study dependence of parameters of the curves. As we know
from the theory of integrable systems, (see \cite{DN}), it is
important to study not only a single curve and its Jacobian, but
rather a whole class of curves and their Jacobians. Following this
line,
 let $B$ denote  a set  parameterizing nonsingular
curves of genus 2 and let $U$ be the analytic bundle over $B$ with a
genus 2 curve as a fibre. We denote by $J(U)$ the associated bundle
over $B$ with a Jacobian of a genus 2 curve as  a fiber. Then,
$J(U)$ has a natural structure of a topological groupoid, see \cite
{BL}.

In the next section, after recalling the definition of the $n$
valued group, we are introducing a notion of $n$-groupoid, which
enables us to consider not only a single structure of a two-valued
group, but also dependance of the parameters of the curve. Thus the
notion of $n$-groupoid is a natural analogue of the notion of
topological groupoid.

Another important reason of using $\sigma$-functions, lies in the
fact that addition relations for them are well adjusted to
two-valued structure, as one can easily see from the genus one case
formula:
\[ \sigma(u+v)\sigma(u-v) = \sigma(u)^2\sigma(v)^2 \big( \wp(v)-\wp(u) \big), \]

The program of construction of $\sigma$-functions of higher genera
had been proposed by Klein. In genus 2 case, the program got rather
advanced development by Baker (see \cite {Bak}), although in his
last review of the program in 1923, Klein noted that program for
genus $g>2$ was still far from being completed. We base our first
approach on recent results, see \cite {BL} and references therein.

Our second  approach is geometric. Following \cite {Don}, \cite
{Kn}, \cite {Re} a geometric approach to addition laws on
hyperelliptic Jacobians has been developed further in \cite {DR2}
(see also \cite {DR3}). It has been  based on the notion of billiard
algebra and it has been connected to integrable billiard systems
within confocal quadrics.

Our first approach is developed in Section \ref{sec:kummersigma},
and the second one in Section \ref{sec:kummergeom}. The Section
\ref{sec:kummersigma} is preceded by the three sections of
preparation. In Section \ref{sec:cp1}, the derivation of the
two-valued group law on $W/\tau$, where $W$ denotes an elliptic
curve as above, is based on $\sigma$ - functions. The rational limit
of the last structure is studied in the next Section
\ref{sec:elementary}. The study of the rational limit of a Kummer
surface and related two-valued group law is performed in Section
\ref{sec:ratkummer}.

In Section \ref{sec:int}, the geometric approach from the Section
\ref{sec:kummergeom}, develops further toward integrable systems,
through the integrable billiards and the billiard algebra and is
motivated by \cite {DR2} .

The final Section \ref{sec:moduli} gives yet another construction of
the two-valued group law on the Kummer variety,  based on the
realization of the Kummer variety in the framework of moduli of
semi-stable bundles. This framework has been developed by Narasimhan
and Ramanan, see \cite {NR}. This last construction of the
two-valued structure is remarkably similar to the historically first
example of topological two-valued group from 1971 (see \cite{BN}),
with a significant difference. The resulting bundles in the 1971
case of \cite {BN} are "virtual", while in the present case the
resulting bundles are effectively realizable.

\medskip

\section{From $n$-valued groups to $n$-groupoids.}\label{sec:def}
\medskip

\subsection{Defining notions and basic examples of multivalued
groups}\label{subsection:defexample}

\medskip

Following \cite{Buc}, we give the definition of an n-valued group on
$X$ as a map:
$$
\aligned
 &m:\, X\times X \rightarrow (X)^n\\
 &m(x,y)=x*y=[z_1,\dots, z_n],
 \endaligned
 $$
where $(X)^n$ denotes the symmetric $n$-th power of $X$ and $z_i$
coordinates therein.

{\it Associativity} is the condition of equality of two $n^2$-sets
$$
\aligned &[x*(y*z)_1,\dots, x*(y*z)_n]\\
&[(x*y)_1*z,\dots, (x*y)_n*z]
\endaligned
$$
for all triplets $(x,y,z)\in X^3$.

An element $e\in X$ is {\it a unit} if
$$
e*x=x*e=[x,\dots,x],
$$
for all $x\in X$.

A map $\inv: X\rightarrow X$ is {\it an inverse} if it satisfies
$$
e\in \inv(x)*x, \quad e\in x*\inv(x),
$$
for all $x\in X$.

Following \cite{Buc}, we say that $m$ defines {\it an $n$-valued
group structure} $(X, m, e, \inv)$ if it is associative, with a unit
and an inverse.

An $n$-valued group $X$ acts on the set $Y$ if there is a mapping
$$
\aligned &\phi:\, X\times Y \rightarrow (Y)^n\\
&\phi (x,y)=x\circ y,
\endaligned
$$
such that the two $n^2$-multisubsets of $Y$
$$
x_1\circ (x_2\circ y) \quad (x_1*x_2)\circ y
$$
are equal for all $x_1, x_2\in X, y\in Y$. It is additionally
required that
$$
e\circ y=[y,\dots, y]
$$
for all $y\in Y$.
\medskip
\begin{example}\label{ex:zplus}[A two-valued group structure on $\mathbb {Z}_+$, \cite{BV}]
Let us consider the set of nonnegative integers $\mathbb Z_+$ and
define a mapping
$$
\aligned &m:\, \mathbb Z_+ \times \mathbb Z_+ \rightarrow (\mathbb
Z_+)^2,\\
&m(x,y)=[x+y,|x-y|].
\endaligned
$$
This mapping provides a structure of a two-valued group on $\mathbb
Z_+$ with the unit $e=0$ and the inverse equal to the identity
$\inv(x)=x$.

In \cite{BV} sequence of two-valued mappings associated with the
Poncelet porism was identified as the algebraic representation of
this 2-valued group. Moreover, the algebraic action of this group on
$\mathbb {C}P^1$ was studied and it was shown that in the
irreducible case all such actions are generated by Euler-Chasles
correspondences.
\end{example}
\medskip
\begin{example}\label{ex:cplus}[2-valued group on $(\mathbb C,+)$]
\medskip
Among the basic examples of multivalued groups, there are $n$-valued
additive group structures on $\mathbb C$. For $n=2$, this is a
two-valued group $p_2$ defined by the relation
\begin{equation}\label{eq:p2}
\aligned &m_2:\, \mathbb C \times \mathbb C \rightarrow (\mathbb
C)^2\\
&\, x *_2 y =[(\sqrt{x}+\sqrt{y})^2, (\sqrt{x}-\sqrt{y})^2]
\endaligned
\end{equation}
\medskip
The product $x *_2 y$ corresponds to the roots in $z$ of the
polynomial equation
$$
p_2(z, x, y)=0,
$$
where
$$
p_2(z, x, y)= (x+y+z)^2-4(xy+yz+zx).
$$

This two-valued group structure has been connected with
degenerations of the Kowalevski top in \cite{Drag2}. Similar
integrable systems were studied by Appel'rot, Delone, Mlodzeevskii
(see \cite{Del}, \cite{App}, \cite {Mlo}, \cite{Gol}.)

As it has been observed in \cite {Drag2}, the general Kowalevski
case is connected with $p_2$ together with its deformation on
$\mathbb {C}P^1$ as a factor of an elliptic curve, see the next
example.

\end{example}
\medskip
\begin{example}\label{ex:eplus}[2-valued group on $S^2=\mathbb {\hat C}$, associated with an elliptic curve]
\medskip
 Suppose
that a cubic  $W$ is given in the standard form
$$
W: t^2=J(s)=4s^3-g_2s-g_3.
$$
Consider the mapping $W \rightarrow S^2 = \mathbb {\hat C} : (s,t)
\mapsto s$, where $\mathbb {\hat C}$ represents a complex line
extended by $\infty$.

The curve $W$ as a cubic curve has the group structure. Together
with its canonical involutive automorphism $\tau: (s,t)\mapsto
(s,-t)$, it defines the standard two-valued group structure of coset
type (see \cite {BR}, \cite{Buc}) on $S^2$ with the unit at infinity
in $S^2$. The product is defined by the formula:
\begin{equation}\label{eq:G2Z2}
[s_1] *_c [s_2] =
\left[\left[-s_1-s_2+\left(\frac{t_1-t_2}{2(s_1-s_2)}\right)^2\right],
\left[-s_1-s_2+\left(\frac{t_1+t_2}{2(s_1-s_2)}\right)^2\right]\right],
\end{equation}
where $t_i=J(s_i), i=1,2$, and
$$[s_i]=\{(s_i,t_i),(s_i,-t_i)\}, \quad s_i=\wp(u_i), t_i=\wp'(u_i),$$
by using addition theorem for the Weierstrass function $\wp(u)$:
$$
\wp(u_1+u_2)=-\wp(u_1)-\wp(u_2) +
\left(\frac{\wp'(u_1)-\wp'(u_2)}{2(\wp(u_1)-\wp(u_2))}\right)^2.
$$

\medskip
The Kowalevski integration procedure  was explained in \cite{Drag2}
as certain deformation of $p_2$ to $(W, \tau)$.

The 2-group structure $(W,\tau)$ is also connected with  the
Poncelet and the Darboux theorem (see \cite{Dar1}, \cite{Drag1},
\cite{DR2}).

In this example, $g_2, g_3$ are parameters of the curve, and they
lead to a rational limit, in the standard limit procedure.

Related structure on $\mathbb {C}P^1$ as an algebraic  mapping
$\mathbb {C}P^1\times \mathbb {C}P^1\rightarrow \mathbb {C}P^2$ (see
\cite {BR}, \cite{Buc}) is presented in Section \ref{sec:cp1}.
\end{example}

\medskip
\begin{example}\label{ex:bundles}

Let $X$ be a topological space. Denote by $W(X)$ the set of all
2-dimensional complex vector bundles $\zeta=\eta +\bar\eta$ over
$X$, where $\eta$ is a linear complex vector bundle over $X$.
 Then the formula

$$\zeta_1\otimes\zeta_2=(\eta_1\otimes\eta_2+\bar\eta_1\otimes\bar\eta_2)+
(\bar\eta_1\otimes\eta_2+\eta_1\otimes\bar\eta_2)$$

 gives the structure of 2-valued group on $W(X)$ defined by the 2-valued
multiplication

$$\zeta_1\star\zeta_2=[(\eta_1\otimes\eta_2+\bar\eta_1\otimes\bar\eta_2),
(\bar\eta_1\otimes\eta_2+\eta_1\otimes\bar\eta_2)].$$

\end{example}

\medskip

\subsection{Topological $n$-groupoids}

\medskip

Following \cite {BL}, let us fix a topological space $Y$. A space
$X$ with a map $p_X: X\rightarrow Y$ is called a space over $Y$, and
the map $p_X$ is called anchor. For $Y$, the anchor $p_Y$ is the
identity map.

A map $f:X_1\rightarrow X_2$ of two spaces over $Y$ is called a map
over $Y$ if $p_{X_2}\circ f=p_{X_1}$. For two spaces $X_1$, $X_2$
over $Y$, their direct product over $Y$ is defined as
$$
X_1\times_Y X_2=\{(x_1, x_2)\in X_1\times X_2\mid
p_{X_1}(x_1)=p_{X_2}(x_2)\}
$$
with the map $p_{X_1\times_Y X_2}(x_1, x_2)=p_{X_1}(x_1)$. Along
this line, one may define a product over $Y$ of $n$ spaces $X_1,
\dots X_n$ over $Y$
$$
X_1\times_Y\dots \times_Y X_n=\{(x_1,\dots, x_n)\in X_1\times
X_2\times \dots \times X_n \mid
p_{X_1}(x_1)=p_{X_2}(x_2)=\dots=p_{X_n}(x_n)\}.
$$
In a special case $X_1=\dots=X_n=X$, we define {\it $n$-th power of
$X$ over $Y$} and denote it as $X^n_Y$.

For a space $X$ over $Y$ we define its {\it $n$-th symmetric power
over $Y$}, denoted as
$$
(X)^n_Y
$$
as the quotient of $X^n_Y$ by the action
of the permutation group.

We define also the diagonal map over $Y$
$$
D:X\rightarrow (X)^n_Y, \quad x\mapsto (x, x, \dots, x).
$$

\begin{definition} A space $X$ with an anchor  $p_X:X\rightarrow Y$
and structural maps over $Y$:
$$
\mu: X\times_Y X\rightarrow (X)^n_Y, \quad inv:X\rightarrow X, \quad
e:Y\rightarrow X
$$
is called \emph{ an $n$-groupoid over $Y$} if the following
conditions are satisfied:
\begin{itemize}
\item [1] For $x_1, x_2, x_3$ such that
$p_X(x_1)=p_X(x_2)=p_X(x_3)$: if
$$
\mu(x_1, x_2)=[z_1,\dots, z_n], \quad \mu(x_2, x_3)=[w_1,\dots,
w_n],
$$
then
$$
[\mu(x_1, w_1), \dots, \mu(x_1, w_n)]=[\mu(z_1, x_3), \dots,
\mu(z_n, x_3)].
$$
\item[2] For every $x\in X$ and $y=p_X(x)\in Y$:
$$
\mu(e(y), x)=\mu(x, e(y))= D(x).
$$
\item[3] For every $x\in X$ and $y=p_X(x)\in Y$:
$$
e(y) \in \mu(x, inv(x)), \quad e(y)\in \mu(inv(x), x).
$$
\end{itemize}
\end{definition}

\medskip
\begin{example} Let $X=\mathbb C\times\mathbb C$, $Y=\mathbb C$ and
an anchor is defined as projection to the second component
$$
p_X:=p_2: X\rightarrow Y, \quad p_X(x, \lambda):=\lambda.
$$
We define $2$-groupoid over $Y$ starting with an operation over $Y$
$$
\mathcal A((x_1, \lambda),(x_2, \lambda))=(x_1+x_2-\lambda x_1x_2,
\lambda).
$$
An involutive automorphism $I$ over $Y$ is defined by
$$
I: X\rightarrow X,\quad {(\bar u,\lambda)}=\left(-\frac
{u}{1-\lambda u}, \lambda\right).
$$
If we denote by $u \bar u = x,\, v\bar v=y$, then the $2$-groupoid
over $Y$ is defined by
$$
\mu ((x, \lambda), (y,\lambda))=[(z_1, \lambda), (z_2,\lambda)]
$$
where $z_i,\, i=1,2$ are solutions of the following quadratic
equation
$$
Z^2-\left(2(x+y)-\lambda^2xy\right)Z+(x-y)^2=0.
$$
Notice that we get $p_2$ for $\lambda=0$. Thus the structure of the
above 2-groupoid is certain deformation of $p_2$.

\end{example}
\medskip
\begin{example}
Similarly, we can consider a two-parameter deformation of elementary
two-valued group $p_2$ and we  get  a structure of 2-groupoid over
$Y_1=\mathbb {C}^2$. We are starting with an operation over $Y_1$ on
$X_1=\mathbb C\times Y_1$, defined by
$$
\mathcal A_1\left((u, \lambda_1, \lambda_2), (v,\lambda_1,
\lambda_2)\right):=\left(\frac{u+v-\lambda_1uv}{1-\lambda_2uv},
\lambda_1, \lambda_2\right).
$$
The involutive automorphism $I_1$, over $Y_1$, is defined by
$$
I_1(u, \lambda_1,\lambda_2)=\left(-\frac{u}{1-\lambda_1
u},\lambda_1, \lambda_2\right).
$$
One can easily deduce corresponding morphism $\mu_1$:
 $$ \mu_1 ((x,
\lambda_1, \lambda_2), (y,\lambda_1, \lambda_2))=[(z_1, \lambda_1,
\lambda_2), (z_2,\lambda_1, \lambda_2)]
$$
where $z_i,\, i=1,2$ are solutions of the following quadratic
equation:
$$
Z^2-\frac{B(x,y, \lambda_1,\lambda_2)}{G(x,y,\lambda_1,
\lambda_2)}Z+\frac{C(x,y,\lambda_1,
\lambda_2)}{G(x,y,\lambda_1,\lambda_2)}=0.
$$
Here we use the notation
$$x=uI_1(u), \quad y=vI_1(v)$$
and
$$
\aligned
 B(x,y, \lambda_1,\lambda_2)&=x^2y^2\lambda_1^2\lambda_2(-\lambda_1^2+\lambda_2)+
 xy(x+y)\lambda_2(2\lambda_2-3\lambda_1)+2(x+y) \\
C(x,y,\lambda_1,\lambda_2)&=(x-y)^2\\
G(x,y,\lambda_1,\lambda_2)&=1+xy\lambda_2(\lambda_1^2-\lambda_2)+xy(x+y)\lambda_1^2\lambda_2^2+
x^2y^2\lambda_2^3(\lambda_2-\lambda_1^2)
\endaligned
$$
\end{example}
\medskip

\section{Structure of two-valued group on  $\mathbb{CP}^1$ and
sigma-functions}\label{sec:cp1}

\medskip

 Let us consider the Weierstrass sigma-function  $\sigma(u) = \sigma(u,g_2,g_3)$, associated with a curve
\[ V=\{ (t,s)\in\mathbb{C}^2\;:\;t^2=4s^3-g_2s-g_3 \}. \]
On the Jacobian  $J_1 = \mathbb{C}^1/\Gamma$ of the curve $V$, the
$\wp$-function $\wp(u)=-\frac{d^2}{du^2}\ln\sigma(u)$ is defined.
The function $\sigma(u)$ is odd, thus, the mapping
\[ \pi \colon J_1 \longrightarrow \mathbb{C}P^1\;:\; \pi(u) =(x_1:x_2), \]
with $x_1=\sigma(u)^2$ and $x_2=\sigma(u)^2\wp(u)$, factorizes as a
composition
\[ \pi \colon J_1 \longrightarrow  J_1/_{\pm} \stackrel{\widehat\pi}{\longrightarrow} \mathbb{C}P^1, \]
where $\widehat\pi$ is a homeomorphism. Let us note that $x_1(u)$
and $x_2(u)$ are entire functions of $u$.

The canonical homeomorphism
\[ \gamma \colon (\mathbb{C}P^1)^2 \longrightarrow \mathbb{C}P^2 \;: \;
\big((x_1:x_2),(y_1:y_2)\big)\longrightarrow
\big(x_1y_1:(x_1y_2+x_2y_1):x_2y_2\big), \] corresponds to the
mapping
\[ \big[ (x_1t_1+x_2t_2),(y_1t_1+y_2t_2) \big] \longrightarrow x_1y_1t_1^2+(x_1y_2+x_2y_1)t_1t_2+x_2y_2t_2^2. \]

\begin{theorem} Multiplication in the two-valued group
\[ m \colon  J_1/_{\pm} \times J_1/_{\pm} \longrightarrow( J_1/_{\pm})^2 \]
is defined  by algebraic mapping
\[  m \colon \mathbb{C}P^1 \times \mathbb{C}P^1  \longrightarrow
\mathbb{C}P^2 \;:\; (x_1:x_2)\ast(y_1:y_2) = (z_1:z_2:z_3), \]
 where $z_1 = (x_1y_2-x_2y_1)^2,\;\, z_2 = 2 \left[ (x_1y_2+x_2y_1)
 \big( x_2y_2-\frac{g_2}{4}x_1y_1 \big)-\frac{g_3}{2}x_1^2y_1^2
 \right],\;\\
 z_3 = \big( x_2y_2+\frac{g_2}{4}x_1y_1 \big)^2 + g_3x_1y_1(x_1y_2+x_2y_1).$
\end{theorem}

\begin{proof}
The composition of mappings
\[  J_1/_{\pm} \times J_1/_{\pm} \stackrel{m}{\longrightarrow}
(J_1/_{\pm})^2 \stackrel{(\widehat\pi)^2}{\longrightarrow}
(\mathbb{C}P^1)^2 \stackrel{\gamma}{\longrightarrow} \mathbb{C}P^2
\]
has the form
\[ [u]\times [v] \longrightarrow \big( [u+v],[u-v] \big) \longrightarrow
\big((x_1:x_2),(y_1:y_2)\big)\longrightarrow (z_1:z_2:z_3), \] where
$[u]=\{ u,-u \},\;\, x_1=\sigma(u+v)^2,\;\,
 x_2=x_1\wp(u+v),\;\,
 y_1=\sigma(u-v)^2,\;\,  y_2=y_1\wp(u-v)$.\\ Thus,
 $$ z_1=\sigma(u+v)^2\sigma(u-v)^2,\quad z_2=z_1\big( \wp(u+v)+\wp(u-v)
 \big),\quad  z_3=z_1\wp(u+v)\wp(u-v).$$

Our goal is to express  $z_1,z_2,z_3$ as functions of
$x_1,x_2,y_1,y_2$. By using the classical addition theorem for
sigma-functions
\[ \sigma(u+v)\sigma(u-v) = \sigma(u)^2\sigma(v)^2 \big( \wp(v)-\wp(u) \big), \]
we get
\[  z_1=(x_1y_2-x_2y_1)^2. \]

From
\[ \ln\sigma(u+v)+\ln\sigma(u-v) = 2\big( \ln\sigma(u)+\ln\sigma(v)
\big)+\ln\big( \wp(v)-\wp(u) \big), \] applying the operator
$\partial=\frac{\partial}{\partial u} + \frac{\partial}{\partial
v}$, we get
\[ 2\zeta(u+v) = 2\big( \zeta(u)+\zeta(v) \big) + \frac{ \wp'(v)-\wp'(u)}{\wp(v)-\wp(u)} \]
where $\zeta(u)=\frac{\partial}{\partial u}\ln\sigma(u)$. By
differentiating previous identity in $u$, one gets
\[ \wp(u+v) = w_1(u,v)-w_2(u,v), \]
where
\begin{equation}\label{eq:w1}
w_1(u,v) = \wp(u)+\frac{1}{2}\frac{\wp''(u)\big(
\wp(v)-\wp(u)\big)+\wp'(u)^2}{\big(\wp(v)-\wp(u)\big)^2}
\end{equation}
\begin{equation}
w_2(u,v) = \frac{1}{2}
\frac{\wp'(u)\wp'(v)}{\big(\wp(v)-\wp(u)\big)^2}.
\end{equation}
Observe, that $w_1(u,-v) = w_1(u,v)$ and  $w_2(u,-v) = -w_2(u,v)$.
Finally we get,
\[ \wp(u+v)+\wp(u-v) = 2w_1(u,v); \quad  \wp(u+v)\wp(u-v) = w_1^2-w_2^2. \]
By using the Weierstrass uniformization of the elliptic curve, we
get
\[ \wp'(u)^2 = 4\wp(u)^2-g_2\wp(u)-g_3, \]
and also,
\[ \wp''(u) = 6\wp(u)^2-\frac{g_2}{2}. \]

According to the formula \ref{eq:w1} we get
\begin{multline*}
2\big(\wp(v)-\wp(u)\big)^2 w_1(u,v) =\\
=2\wp(u)\big(\wp(v)-\wp(u)\big)^2 + \left( 6\wp(u)_2-\frac{g_2}{2}
\right) \big(\wp(v)-\wp(u)\big) +4\wp(u)^3 - g_2\wp(u)-g_3 =\\
= 2\wp(u)\wp(v)\big(\wp(v)+\wp(u)\big) -
\frac{g_2}{2}\big(\wp(v)+\wp(u)\big)-g_3.
\end{multline*}
Thus,
\begin{align*}
2z_1^2w_1(u,v) &=
2\sigma^4(u)\sigma^4(v)\big(\wp(v)-\wp(u)\big)^2w_1(u,v) =\\
&= \sigma^4(u)\sigma^4(v)\left[
\big(\wp(v)+\wp(u)\big)\Big(2\wp(v)\wp(u)-\frac{g_2}{2}\Big)-g_3
\right].
\end{align*}
 As a consequence, we get
\[ z_2 = 2\left[
(x_1y_2+x_2y_1)\Big(x_2y_2-\frac{g_2}{4}\,x_1y_1\Big)-\frac{g_3}{2}\,x_1^2y_1^2
\right].\]
 Further, we get
\begin{multline*}
(\wp(v)-\wp(u)\big)^4(w_1^2-w_2^2) = \\
= \left[
\big(\wp(v)+\wp(u)\big)\Big(2\wp(v)\wp(u)-\frac{g_2}{4}\Big)-\frac{g_3}{2}
\right]^2 - \left[ 2\wp(u)^3-\frac{g_2}{2}\wp(u)-\frac{g_3}{2}
\right]\left[ 2\wp(v)-\frac{g_2}{2}\wp(v)-\frac{g_3}{2} \right].
\end{multline*}
From
\begin{multline*}
 z_1 z_3 = \sigma^8(u)\sigma^8(v)\big(\wp(v)-\wp(u)\big)^4(w_1^2-w_2^2) = \\
\left[
(x_1y_2+x_2y_1)\Big(x_2y_2-\frac{g_2}{4}\,x_1y_1\Big)-\frac{g_3}{2}\,x_1^2y_1^2
\right]^2- \\
-x_1y_1\left[ 2x_2^3 - \frac{g_2}{2}\,x_2x_1^2 - \frac{g_3}{2}
\,x_1^3 \right] \left[ 2y_2^3 - \frac{g_2}{2}\,y_2y_1^2 -
\frac{g_3}{2} \,y_1^3 \right],
\end{multline*}
we get
\[  z_1 z_3 = \left[
\Big(x_2y_2+\frac{g_2}{4}\,x_1y_1\Big)^2+g_3x_1y_1(x_1y_2+x_2y_1)
\right]. \]
\end{proof}

\begin{corollary}
The multiplication law of the two-valued group $\mathbb{C}P^1$ is
defined in homogeneous coordinates by the formula
\begin{multline*}
(x_1t_1+x_2t_2)\ast(y_1t_1+y_2t_2) =\\
=(x_1y_2-x_2y_1)^2t_1^2 + 2\left[ (x_1y_2+x_2y_1)
\Big(x_2y_2-\frac{g_2}{4}\,x_1y_1\Big)-\frac{g_3}{2}\,x_1^2y_1^2
\right]t_1t_2 +\\
+\left[\Big(x_2y_2+\frac{g_2}{4}\,x_1y_1\Big)^2 +
g_3x_1y_1(x_1y_2+x_2y_1) \right]t_2^2.
\end{multline*}
In the rational limit, defined with  $g_2=g_3=0$, we get
\begin{multline*}
(x_1t_1+x_2t_2)\ast(y_1t_1+y_2t_2) =\\
=(x_1y_2-x_2y_1)^2t_1^2 + 2(x_1y_2+x_2y_1)x_2y_2t_1t_2 +
x_2^2y_2^2t_2^2 =\\
= \left[(x_1y_2+x_2y_1)t_1 + x_2y_2t_2 \right]^2
-4x_1x_2y_1y_2t_1^2.
\end{multline*}
\end{corollary}

\medskip

\medskip

\section{Elementary two-valued group on $(\mathbb C^2,
+)$}\label{sec:elementary}

\medskip

Let us consider $\mathbb C^2$ together with its subgroup $\mathbb
Z^2$. By the standard coset construction (see \cite{Buc}),
$$X=\mathbb C^2/\pm$$
is equipped by two-valued group structure:
$$
\mu:\, X\times X\rightarrow (X)^2
$$
given in coordinates $u=(u_1,u_2), v=(v_1,v_2)\in \mathbb C^2$ by
the formula
$$
\mu ([u,-u],[v,-v])=([u+v, - (u+v)], [(u-v),-(u-v)]).
$$
\medskip

Now, we pass to the embedding
$$
\pi:\, X\rightarrow \mathbb C^3, \quad \pi ([(u_1,
u_2),-(u_1,u_2)])=(x_1,x_2,x_3)
$$
given by the formula
$$
x_1=u_1^2, \quad x_2=u_1u_2 \quad x_3=u_2^2.
$$
By checking, one sees that the inverse image of a point is a two
element set: $\pi^{-1}(x_1,x_2,x_3)=[(u_1,u_2),-(u_1,u_2)]$, where
$$
u_1=\pm \sqrt{x_1}, \quad u_2=\pm \sqrt{x_3}.
$$
The image of $\pi$ is a quadric $Q$ in $\mathbb C^3$ defined by the
equation
$$
Q: x_1x_3=x_2^2.
$$
On $Q$, the multiplication $\mu$ can be rewritten according to the
formulae
$$
\mu ([u,-u],[v,-v])=((X_2, X_4, X_6), (Y_2,Y_4,Y_6)),
$$
where
$$
[u,-u]=\pi ^{-1}(x_1,x_2,x_3),\quad [v,-v]=\pi ^{-1}(y_1,y_2,y_3).
$$
Indices $2, 4, 6$ for variables $X$ and $Y$ are chosen to fit with a
graduation in a sequel. Then
$$
X_2=(\sqrt{x_1}+\sqrt {y_1})^2, \quad Y_2=(\sqrt{x_1}-\sqrt
{y_1})^2,
$$
and $X_2, Y_2$ are the roots of the equation
\begin{equation}\label{eq:eg1}
Z_1^2-2(x_1+y_1)Z_1+(x_1-y_1)^2=0.
\end{equation}
For $X_4$ and $Y_4$ we have
$$
\aligned X_4&=(\sqrt{x_1}+\sqrt {y_1})(\sqrt{x_3}+\sqrt {y_3})\\
Y_4&=(\sqrt{x_1}-\sqrt {y_1})(\sqrt{x_3}-\sqrt {y_3}).
\endaligned
$$
Thus,
$$\aligned
X_4+Y_4&=2(x_2+y_2)\\
X_4Y_4&=(x_1-y_1)(x_3-y_3),
\endaligned
$$
and $X_4, Y_4$ are the solutions of the equation:
\begin{equation}\label{eq:eg2}
Z_2^2-2(x_2+y_2)Z_2+(x_1-y_1)(x_3-y_3)=0.
\end{equation}
The last equation can be written in the form
$$
(Z_2-(x_2+y_2))^2=y_1x_3+x_1y_3 + 2x_2y_2.
$$
Finally, for $X_6, Y_6$ the situation is analogue to $X_2, Y_2$:
$$
X_6=(\sqrt{x_3}+\sqrt {y_3})^2, \quad Y_6=(\sqrt{x_3}-\sqrt
{y_3})^2,
$$
$X_6, Y_6$ are the roots of the equation
\begin{equation}\label{eq:eg3}
Z_3^2-2(x_3+y_3)Z_3+(x_3-y_3)^2=0.
\end{equation}

Three quadratic equations (\ref{eq:eg1}, \ref{eq:eg2}, \ref{eq:eg3})
determine four possible pairs of triplets of solutions. But, the
constraints
$$
X_2X_4=X_6^2, \quad Y_2Y_4=Y_6^2,
$$
select a unique pair of triplets $[(X_2,X_4,X_6), (Y_2,Y_4,Y_6)]$,
which is in accordance with the theory of two-valued groups.

\medskip

\section{Rational limit of a Kummer surface and two-valued addition
law}\label{sec:ratkummer}
\medskip

We start with a genus two curve $V$ given by affine equation
$$
V=\{(x,y)\in \mathbb C^2:\, y^2=x^5+\lambda_4x^3+\lambda_6 x^2+
\lambda_8x+\lambda_{10}\}.
$$
Corresponding Jacobian $\Jac(V)$ as two-dimensional complex torus is
a factor of $\mathbb C^2$ with a lattice $\Gamma$. The lattice
$\Gamma$ is determined by the vector $(\lambda_4, \lambda_6,
\lambda_8, \lambda_{10})$. Note that the indices  $(4, 6, 8, 10)$ of
$\lambda$  are chosen to fit with the graduation in sequel. The
Kummer surface $K$ is the factor of the Jacobian by the group of
automorphisms  of order $2$:
$$
K=\Jac(V)/\pm.
$$
Locally, in a vicinity of $0=(0,0)$ the Kummer surface $K$ is
isomorphic to $\mathbb C^2/\pm$.

Moreover, all constructions allow the rational
limit:$\lambda\rightarrow 0$.

Starting from the vector $(\lambda_4, \lambda_6, \lambda_8,
\lambda_{10})$ and a vector $u=(u_1,u_3)$, a function $\sigma (u,
\lambda)$ is constructed as an entire function in $u$ and $\lambda$,
such that all coefficients $c_{ij}(\lambda)$ in expansion in
$u_1^iu_3^j$ are {\it polynomials} in $\lambda$.

\medskip

\subsection{Two-dimensional Weierstrass functions}

\medskip

Following \cite{BEL}, we introduce
$$
\aligned \zeta_k&=\frac {\partial}{\partial u_k}\ln \sigma (u,
\lambda),
\qquad k=1,3\\
\wp_{kl}&=-\frac {\partial^2}{\partial u_k\partial u_l}\ln \sigma
(u, \lambda), \qquad k,l=1,3.
\endaligned
$$
The last, {\it Weierstrass functions} are functions on the Jacobian
$\Jac(V)$.

We are going to use {\bf the standard} sigma function, the one
defined by  {\it arf-invariants} $(\ell,\ell')$ to be equal to
$$
\ell=(1,1), \quad \ell' =(2,1).
$$
The standard sigma function is odd in $u$,
$$
\sigma(u, \lambda)=-\sigma(-u, \lambda),
$$
while the $\wp$- functions are {\it even}:
$$
\wp_{k,l}(u, \lambda)=\wp_{k,l}(-u,\lambda).
$$
Thus,  the $\wp$- functions  $\wp_{k,l}(u, \lambda)$ generate
well-defined functions on the Kummer surface $K$.

There is an embedding with each fixed $\lambda$:
$$
\aligned \pi_{\lambda}:\, &K\rightarrow \mathbb {C}P^3\\
&[u,-u]\mapsto [\sigma^2(u), \sigma^2(u)\wp_{11}(u),
\sigma^2(u)\wp_{13}(u), \sigma^2(u)\wp_{33}(u)].
\endaligned
$$
Note that $\sigma^2(u)\wp_{kl}(u)$ are entire functions. The
embeddings $\pi_\lambda$ serve to describe addition law on $K$ in
terms of coordinates on $\mathbb {C}P^3$.

One can easily compute the limit of the functions when $\lambda$
tends to zero:
$$
\aligned &\lim_{\lambda\to
0}\sigma(u,\lambda)=\sigma_0(u)=u_3-\frac{1}{3}u_1^3\\
&\zeta_1^{(0)}=\frac {\partial}{\partial u_1}\ln \sigma_0 (u)=-\frac{u_1^2}{\sigma_0(u)}\\
&\zeta_3^{(0)}=\frac {\partial}{\partial u_3}\ln \sigma_0 (u)=-\frac{1}{\sigma_0(u)}\\
&\wp_{11}^{(0)}=-\frac {\partial^2}{\partial u_1^2}\ln \sigma_0(u)=\frac{2u_1\sigma_0(u)+u_1^4}{\sigma_0^2(u)}\\
&\wp_{13}^{(0)}=-\frac {\partial}{\partial u_1}\zeta_3^{(0)} (u)== -\frac{u_1^2}{\sigma_0(u)^2}\\
&\wp_{33}^{(0)}=-\frac {\partial}{\partial u_3}\zeta_3^{(0)}
(u)=\frac{1}{\sigma_0^2(u)}.
\endaligned
$$
Thus,
$$
\sigma_0^2\wp_{13}^{(0)}=-u_1^2,\quad \sigma_0^2\wp_{33}^{(0)}=1.
$$
\medskip
\subsection{Rational limit embedding and two-valued group law}
\medskip
Now, we are going to construct a new two-valued group law on
$\mathbb C^2/\pm$ associated with an embedding $\pi_K$ induced by
the rational limit of a Kummer surface:
$$
\aligned &\pi_K:\, \mathbb C^2/\pm \rightarrow \mathbb C^3\\
&\pi_K([u,-u])=\left((u_3-\frac{1}{3}u_1^3)^2,
2u_1u_3+\frac{1}{3}u_1^4,-u_1^2\right),
\endaligned
$$
where
$$
u=(u_1,u_3).
$$
One checks that an inverse image is a two element set if nonempty:
$$
\aligned x_1&=\left(u_3-\frac{1}{3}u_1^3\right)^2\\
x_2&=2u_1u_3+\frac{1}{3}u_1^4\\
x_3&=-u_1^2.
\endaligned
$$

 The multiplication $\mu$ given by the formula
$$
\aligned
 \mu ([u,-u], [v, -v])=&[(u_1+v_1), -(u_1+v_1),
(u_3+v_3),-(u_3+v_3),\\
&(u_1-v_1), -(u_1-v_1), (u_3-v_3),-(u_3-v_3)]
\endaligned
$$
after composition with the embedding $(\pi_K)^2$ leads to the
formulae
$$
\aligned \hat X_6&=\left((u_3+v_3)-\frac{1}{3}(u_1+v_1)^3\right)^2\\
\hat X_4&=2(u_1+v_1)(u_3+v_3)+\frac{1}{3}(u_1+v_1)^4\\
\hat X_2&=-(u_1+v_1)^2
\endaligned
$$
and
$$
\aligned \hat Y_6&=\left((u_3-v_3)-\frac{1}{3}(u_1-v_1)^3\right)^2\\
\hat Y_4&=2(u_1-v_1)(u_3-v_3)+\frac{1}{3}(u_1-v_1)^4\\
\hat Y_2&=-(u_1-v_1)^2.
\endaligned
$$
The last formulae lead to the following change of variables:
$$
\aligned \hat X_2&=-X_2\\
\hat X_4 &= 2 X_4+\frac{1}{3}X_2^2\\
\hat X_6 &= X_6-\frac{2}{3}X_2X_4+\frac{1}{9}X_2^3.
\endaligned
$$
This is an {\it algebraic} change of variables and the inverse
change is given by the formulae:
$$
\aligned X_2&=-\hat X_2\\
X_4&=\frac{1}{2}\left(\hat X_4-\frac{1}{3}\hat X_2^2\right)\\
X_6&=\hat X_6-\frac{1}{3}\hat X_2\hat X_4-\frac{2}{9}X_2^3
\endaligned
$$
\medskip
By applying the last algebraic change of variables on the equation
of the quadric $Q: X_2X_6=X_4^2$ we get the equation of the rational
limit of the Kummer surface.

\medskip
\begin{proposition}\label{prop:ratkum}
 The rational limit of the Kummer surface is
given by the surface in $\mathbb C^3$ by the equation
\begin{equation}\label{eq:ratkum}
-9\hat X_4^2-36\hat X_2\hat X_6+12\hat X_2^2\hat X_4 + 7\hat
X_2^4=0.
\end{equation}
\end{proposition}

\subsection{Rational Kummer two-valued group}
\medskip

In the previous section, we have constructed a two-valued group law
in coordinates $(X_2, X_4, X_6)$ - the elementary two-valued group.
Now, using the algebraic change of variables, we are going to
construct a new two-valued group law, and we will call it {\it the
rational Kummer two-valued group}.

First, we consider $(\hat X_2, \hat Y_2)$. We have
$$
\aligned \hat X_2&=-x_3-2u_1v_1-y_3\\
\hat Y_2&=-x_3+2u_1v_1-y_3
\endaligned
$$
therefore
$$
\aligned \hat X_2+\hat Y_2&=-2(x_3+y_3)\\
\hat X_2\hat Y_2&=(x_3-y_3)^2.
\endaligned $$

Thus we see that the pair $(\hat X_2, \hat Y_2)$ is the solution of
the quadratic equation
\begin{equation}\label{eq:rk1}
\mathcal Z^2+2(x_3+y_3)\mathcal Z +(x_3-y_3)^2=0.
\end{equation}
\medskip

Now, we pass to the pair $(\hat X_4, \hat Y_4)$. They can be
represented in the form
$$
\aligned \hat X_4&=\hat X_4^+ +\hat X_4^-\\
\hat Y_4&=\hat X_4^+ -\hat X_4^-
\endaligned
$$
where
$$
\aligned \hat
X_4^+&=2(u_1u_3+v_1v_3)+\frac{1}{3}(u_1^4+6u_1^2v_1^2+v_1^4)\\
\hat X_4^-&=2(v_1u_3+u_1v_3)+\frac{4}{3}u_1v_1(u_1^2+v_1^2).
\endaligned
$$
Then, we have
$$
\aligned \hat X_4 +\hat Y_4&=2\hat X_4^+\\
\hat X_4\hat Y_4&=(\hat X_4^+)^2-(\hat X_4^-)^2.
\endaligned
$$
Thus,  $(\hat X_4, \hat Y_4)$ are the roots of the quadratic
equation
\begin{equation}\label{eq:rk2}
\mathcal Z^2-2\hat X_4^+\mathcal Z+(\hat X_4^+)^2-(\hat X_4^-)^2=0.
\end{equation}
The last equation is equivalent to
$$
(\mathcal Z-(\hat X_4^+))^2=(\hat X_4^-)^2,
$$
where $\hat X_4^+, \hat X_4^-$ can be rewritten in the form
$$
\aligned \hat X_4^+&=2(u_1u_3+v_1v_3)+\frac{1}{3}(u_1^4+6u_1^2v_1^2+v_1^4) \\
\hat X_4^-&=2(v_1u_3+u_1v_3)+\frac{4}{3}u_1v_1(u_1^2+v_1^2) .
\endaligned
$$
We pass to the last pair $\hat X_6, \hat Y_6$:
$$\aligned
\hat X_6&=\left((u_3+v_3)-\frac{1}{3}(u_1+v_1)^3\right)^2\\
\hat Y_6&=\left((u_3-v_3)-\frac{1}{3}(u_1-v_1)^3\right)^2.
\endaligned
$$
One can easily calculate
$$
\aligned \hat X_6+\hat
Y_6=&2[u_3^2+v_3^2-\frac{2}{3}(u_1^3u_3+3u_1u_3v_1^2+3u_1^2v_1v_3+v_1^4)+\\
&+\frac{1}{9}(u_1^6+15u_1^4v_1^2+15u_1^2v_1^4+v_1^6)]\\
\hat X_6 \hat
Y_6=&\left((u_3^2-v_3^2)+\frac{1}{9}(u_1^2-v_1^2)^3-\frac{2}{3}(u_3u_1^3+3u_1v_1^2u_3^2-3u_1^2v_1v_3-v_1^3v_3)\right)^2,
\endaligned
$$
or, in the old coordinates
$$
\aligned \hat X_6 +\hat
Y_6&=2\left(x_3+y_3-\frac{2}{3}(x_3x_2+3x_2y_1+3x_1y_2+y_1^2)+\frac{1}{9}(x_1^3+15x_1^2y_1+15x_1y_1^2+y_1^3)\right)\\
\hat X_6 \hat Y_6&=\left((x_3-y_3)+\frac
{1}{9}(x_1-y_1)^2-\frac{2}{3}(x_2x_1+3x_2y_1 - 3 x_1y_2
-y_1y_2)\right)^2.
\endaligned
$$
In the new coordinates one may rewrite
$$
\aligned B_3:=&\hat X_6+\hat Y_6\\
C_3:=&\hat X_6\hat Y_6
\endaligned
$$
and to get finally
$$
\aligned
B_3=&2((\hat x_6-\frac{1}{3}\hat x_4\hat
x_2-\frac{2}{9}\hat x_2^3)+(\hat y_6-\frac{1}{3}\hat y_4\hat
y_2-\frac{2}{9}\hat y_2^3)\\
&-\frac{1}{3}(-\hat x_2\hat x_4+\frac{1}{3}\hat x_2^3-3\hat x_4\hat
y_2+\hat x_2^2\hat y_2-3\hat x_2\hat y_4+\hat x_2 \hat y_2^2+2\hat
y_2^2)\\
&+\frac{1}{9}(-\hat x_2^3+15\hat x_2^2\hat y_2+15\hat x_2 \hat
y_2^2-\hat y_2^3))
\endaligned
$$
and
$$
\aligned C_3=&[(\hat x_6-\frac{1}{3}\hat x_4\hat x_2-\frac{2}{9}\hat
x_2^3-\hat y_6+\frac{1}{3}\hat y_4\hat
y_2+\frac{2}{9}\hat y_2^3)+\frac{1}{9}(\hat y_2-\hat x_2)^3\\
&-\frac{1}{9}(\frac{1}{3}\hat x_2^2-\hat x_2\hat x_4+\hat x_2^2\hat
y_2-3\hat x_4\hat y_2+3\hat y_4\hat x_2-\hat y_2^2\hat x_2\\
&+\hat y_2\hat y_4-\frac{1}{3}\hat y_2^3)]^2.
\endaligned
$$
Thus, we may conclude that the pair $(\hat X_6, \hat Y_6)$ is
determined as the roots of the quadratic equation
\begin{equation}\label{eq:rk3}
\mathcal Z^2-B_3\mathcal Z + C_3=0,
\end{equation}
where $B_3, C_3$ are functions of the coordinates $(\hat x_2, \hat
x_4, \hat x_6, \hat y_2, \hat y_4, \hat y_6)$ given above.

\

\

\section{Two-valued group structures on Kummer varieties and sigma-functions}
\label{sec:kummersigma}
\medskip

We start with the sigma-function $\sigma(u) =
\sigma(u,\lambda)$,where
 $u^\top = (u_1,u_3),\; \lambda^\top = (\lambda_4,\lambda_6,\lambda_8,\lambda_{10})$, associated with a curve
\[ V=\{ (t,s)\in\mathbb{C}^2\;:\;t^2=s^5+\lambda_4s^3+\lambda_6s^2+\lambda_8s+\lambda_{10} \}. \]
We assume that  $\lambda\in \mathbb{C}^4$ is a non-discriminant
point of the curve $V$ and  $u\in \mathbb{C}^2$, where $du_1 =
\frac{sds}{t},\;du_3 = \frac{ds}{t}$. Indexation of the coordinates
of the vector of the  parameter $\lambda$ is chosen according to the
graduation $\deg s = -2,\; \deg t = -5,\; \deg \lambda_{2i} = -2i$.
Moreover, $\deg u_1 = 1,\; \deg u_3 = 3$, and the sigma-function
\[ \sigma(u,\lambda) = u_3 - {1\over 3}\,u_1^3 + {1\over 6}\,\lambda_6u_3^3 + (u^5) \]
is an entire and homogeneous of degree 3 in $u$ and $\lambda$.
Recurrent description of the series for $\sigma(u,\lambda)$ is given
in \cite {BL}.

We consider the Jacobian  $J_2 = \Jac(V) = \mathbb{C}^2/\Gamma_2$ of
the curve  $V$ and a vector-function is defined by the formulae
\[ \wp(u)=\big(\wp_{33}(u),\wp_{13}(u),\wp_{11}(u)\big), \]
where $\wp_{kl}(u) = -\frac{\partial^2}{\partial u_k\partial
u_l}\ln\sigma(u),\;k,l = 1,3$. Since the function $\sigma(u)$ is
odd, we get a mapping
\[ i \colon J_2 \longrightarrow \mathbb{C}P^3\;:\; i(u) =(x_0:x_2:x_4:x_6), \]
where $ x_0=\sigma(u)^2\wp_{33}(u),\; x_2=\sigma(u)^2\wp_{13}(u), \;
x_4=\sigma(u)^2\wp_{11}(u), \; x_6=\sigma(u)^2$. The last mapping
factorizes through the Kummer variety $K$ together with an embedding
\[ \widehat i \colon K=\left( J_2/_{\pm} \right) \longrightarrow  \mathbb{C}P^3. \]
Let us note, that the embedding  $\widehat i$ is defined with entire
homogeneous functions $x_{2k}(u),\; \deg\, x_{2k}=2k,\;
k=0,\ldots,3$.

We are going to use a ramified covering
\[ \gamma \colon \big(\mathbb{C}P^3\big)^2 \longrightarrow \mathbb{C}P^6, \]
defined by the relation
\[ \left[ p(x,t), p(y,t) \right] \longrightarrow p(x,t)p(y,t),  \]
where $p(x,t) = x_0t_3^3+x_2t_3^2t_1+x_4t_3t_1^2+x_6t_1^3$. By
putting  $\deg\, t_k=k$, we get $p(x,t)$ as a homogeneous polynomial
of degree 9.

\begin{theorem} Multiplication $m$ in the two-valued group on $K$
\[ m \colon  K \times K \longrightarrow(K)^2\; : \; [u]\ast[u] = \big( [u+v],[u-v] \big) \]
is defined through algebraic mapping
\[  \mu \colon \mathbb{C}P^3 \times \mathbb{C}P^3  \longrightarrow
\mathbb{C}P^6 \;:\;  \] and it is uniquely defined by the commuting
condition of the following diagram
$$ \xymatrix{
K \times K \ar[r]^m \ar[d]^{\widehat i \times \widehat i}&(K)^2
\ar[r]^{(\widehat i)^2} &
\big(\mathbb{C}P^3\big)^2\ar[d]^{\gamma}\\
\mathbb{C}P^3 \times \mathbb{C}P^3 \ar[rr]^{\mu}&&  \mathbb{C}P^6
 } $$
\end{theorem}

\begin{proof}
Set
\[ X(u) = \big( \wp_{33}(u),\wp_{31}(u),\wp_{11}(u),1 \big)\quad
\text{è} \quad \mathcal{X}(u) = \sigma(u)^2 X(u). \] Consider the
canonical projection
\[ \pi \colon \mathbb{C}^7\backslash 0 \longrightarrow \mathbb{C}P^6\;:\;
\pi(z) = (z_0:z_2:\ldots:z_{12}) = [z]. \] According to the
construction, we have
\[ \widehat i\, [u] = \left[\mathcal{X}(u)\right], \quad
\gamma(\,\widehat i\,)^2\big([u]\ast[v]\big) =
\gamma\big(\left[\mathcal{X}(u+v),\mathcal{X}(u-v)\right]\big). \]
Thus, we have to show that each coordinate
$z_{2k}(u,v),\;k=0,\ldots,6$, of the point
$\gamma\big(\left[\mathcal{X}(u+v),\mathcal{X}(u-v)\right]\big)$ is
a polynomial of the coordinates $x_{2i},\;y_{2i},\; i=0,\ldots,3$ of
the points  $\mathcal{X}(u)$ è $\mathcal{X}(v)$.

Genus two sigma-function $\sigma(u)$ satisfies the following
addition theorem (see \cite{BEL}, \cite{BL}):
\[ \sigma(u+v)\sigma(u-v) = \mathcal{X}(u)^\top\mathcal{J}\,\mathcal{X}(v). \]
where $\mathcal{J} =
\begin{pmatrix}
0&-\mathcal{E}\\
\mathcal{E}&0
\end{pmatrix}$
and $\mathcal{E} =
\begin{pmatrix}
0&1\\
1&0
\end{pmatrix}$. We have,
\[ z_{12} = z_{12}(u,v) = \big(\mathcal{X}(u)^\top\mathcal{J}\,\mathcal{X}(v)\big)^2. \]
Thus, we get that in the mapping $\mu(x,y) = [z] =
(z_0:z_2:\ldots:z_{12}) \in \mathbb{C}P^6$, the coordinate $z_{12}$
is defined by the formula $z_{12} = (x^\top \mathcal{J} y)^2$, where
$(x,y)\in \mathbb{C}P^3 \times \mathbb{C}P^3 $.

Analogous result for the rest of the coordinates is based on deep
facts about algebraic generators of the ring generated by
logarithmic derivatives of order 2 and higher of the sigma-function
$\sigma(u)$ (see \cite{BEL}).

Set
\[ M(u,v) = X(u)^\top\mathcal{J}\,X(v). \]
We have
\[ \ln\sigma(u+v)+\ln\sigma(u-v) = 2\big(\ln\sigma(u)+\ln\sigma(v)\big)+\ln M(u,v). \]
Apply the operators $\partial_k = \frac{\partial}{\partial
u_k}+\frac{\partial}{\partial v_k}$, for $k=1$ and 3. We get
\[ 2\zeta_k(u+v) = 2\big(\zeta_k(u)+\zeta_k(v)\big) +\frac{1}{M(u,v)}
\left( X_k(u)^\top\mathcal{J}\,X(v) + X(u)^\top\mathcal{J}\,X_k(v)
\right), \] where $\zeta_k(u) = \frac{\partial}{\partial u_k}
\ln\sigma(u)$ è $X_k(u) = \frac{\partial}{\partial u_k}X(u)$. Now,
we apply the operator $\frac{\partial}{\partial u_l}$, and we get
\[ \wp_{kl}(u+v) = - \frac{\partial}{\partial u_l}\,\zeta_k(u+v)=
\varphi_{kl}(u,v) - \psi_{kl}(u,v), \] where
\begin{align*}
\varphi_{kl}(u,v) &= \wp_{kl}(u)-\frac{1}{2M(u,v)^2}\left\{
\big(X_{kl}(u)^\top\mathcal{J}\,X(v)\big)M(u,v) - X_k(u)^\top
B(v)X_l(u) \right\},\\
\psi_{kl}(u,v) &= \frac{1}{2M(u,v)^2} \big(X_{l}(u)^\top C(u,v)
X_k(v)\big)
\end{align*}
and $B(v) = \mathcal{J}\,X(v)X(v)^\top\mathcal{J}^\top,\quad C(u,v)
= \big(M(u,v) - \mathcal{J}\,X(v)X(u)^\top\big)\mathcal{J}$.

Note $B(-v) = B(v)$ è $C(-u,v) = C(u,-v) = C(u,v)$. Thus,
$\varphi_{kl}(u,-v) = \varphi_{kl}(u,v)$ è $\psi_{kl}(u,-v) =
-\psi_{kl}(u,v)$. We get
\[ \mathcal{X}(u+v) = \sigma(u+v)^2(X_1-X_2)\quad \text{è} \quad
\mathcal{X}(u-v) = \sigma(u-v)^2(X_1+X_2), \] where $X_1 =
(\varphi_{33},\varphi_{13},\varphi_{11},1)$ è $X_2 =
(\psi_{33},\psi_{13},\psi_{11},0)$.
 We obtain
\[ p\big( \mathcal{X}(u+v),t \big)p\big( \mathcal{X}(u-v),t \big)
= z_1\cdot\big( p(X_1,t)^2 - p(X_2,t)^2 \big) \]
 where
\begin{align*}
p(X_1,t) &= \varphi_{33} t_3^3 + \varphi_{13} t_3^2t_1 +
\varphi_{11} t_3t_1^2 + t_1^3,\\
p(X_2,t) &= t_3(\psi_{33} t_3^2 + \psi_{13} t_3t_1 + \psi_{11}
t_1^2).
\end{align*}
From the above formulae for $\varphi_{kl}(u,v)$ è
$\psi_{kl}(u,v)\psi_{pq}(u,v)$, it immediately follows that these
functions are polynomials of
$\wp_{ij}(u),\; \wp_{ijk}(u)\wp_{i'j'k'}(u),\; \wp_{ijpq}(u)$ è \\
$\wp_{ij}(v),\; \wp_{ijk}(v)\wp_{i'j'k'}(v),\; \wp_{ijpq}(v)$. The
functions $\wp_{ijk}(u)\wp_{i'j'k'}(u)$ è $\wp_{ijpq}(u)$, where
$i,j,k,i',j',k',p,q$ take values 1 or 3 independently, are
polynomials of $\wp_{ij}(u)$ (see \cite{BEL}). Consequently, using
 $(x,y)|_{\widehat i \times \widehat i} = \big(
\mathcal{X}(u),\mathcal{X}(v) \big)$, we get coordinates $z_{2k}$ of
the  vector $\mu(x,y)$ as polynomials of the coordinates of the
vectors $x$ and $y$.
\end{proof}

\medskip

\section{Homomorphism of rings of functions, induced by Abel mapping in genus 2}
\label{sec:hom}
\medskip

Let $V = \left\{(s,\mu) \in \mathbb{C}^2\;:\; \mu^2 =
s^5+\lambda_4s^3+\lambda_6s^2 +\lambda_8s+\lambda_{10}\right\}$
denotes a hyperelliptic curve.

\begin{theorem} \label{T-1} The Abel mapping
\[ \mathcal{A} \colon (V)^2 \longrightarrow \Jac V \]
induces a homomorphism of rings of functions
\[ \mathcal{A}^*  \colon \mathcal{F}(\Jac V) \longrightarrow \mathcal{F}\big((V)^2\big) \]
such that
\[ \wp_{11}(u) = s_1+s_2,\quad \wp_{13}(u) = -s_1s_2,\quad \wp_{33}(u)
= \frac{F(s_1,s_2)-2\mu_1\mu_2}{(s_1-s_2)^2}, \]
 where
 $$F(s_1,s_2) =
2\lambda_{10}+\lambda_8(s_1+s_2)+s_1s_2\big(2\lambda_6+\lambda_4(s_1+s_2)\big),
$$
\[ \wp_{111}(u) = 2\frac{\mu_1-\mu_2}{s_1-s_2},\quad
\wp_{113}(u) = 2\frac{s_1\mu_2-s_2\mu_1}{s_1-s_2} \]
\[ \wp_{331}(u) = -2\frac{s_1^2\mu_2-s_2^2\mu_1}{s_1-s_2},\quad
\wp_{333}(u) = 2\frac{\psi(s_1,s_2)\mu_2 -
\psi(s_2,s_1)\mu_1}{(s_1-s_2)^3}, \]

and
$$\psi(s_1,s_2) =
4\lambda_{10}+\lambda_8(3s_1+s_2)+2\lambda_6s_1(s_1+s_2)+
\lambda_4s_1^2(s_1+3s_2)+s_1^3s_2(3s_1+s_2).$$
\end{theorem}

The proof of the Theorem can be found in \cite{BEL}. Let us note
that we use the indexation here different from \cite{BEL}, in
correspondence with the graduation:
$$\deg\, s = -2\,;\; \deg\, \mu = -5\,;\; \deg\, \lambda_{2k} = -2k,\;
k=2,3,4,5\,;\; \deg\, u_i = i,\; i=1,3.$$ This provides additional
opportunity to check the formulae. Observe that $\deg\, \wp_{kl}(u)
= -(k+l),\; \deg\, \wp_{klp}(u) = -(k+l+p)$.

Any Abelian  function on  $\Jac V$ represents a linear function of
$\wp_{111}(u)$ with coefficients which are rational functions of
$\wp_{11}(u)$ and $\wp_{13}(u)$. On the other hand, from the theory
of polysymmetric functions (see \cite{GKZ} and \cite {BR}), it is
known that the field of rational functions on $(\mathbb{C}^2)^2$ in
the coordinates $\left[( s_1,\mu_1), (s_2,\mu_2) \right]$ is
generated by polysymmetric functions
\[ e_{10} = s_1+s_2,\;  e_{20} = s_1s_2,\; e_{01} = \mu_1+\mu_2,\;
e_{02} = \mu_1\mu_2,\;  e_{11} = s_1\mu_2+s_2\mu_1, \] which are
related by unique  (for $(\mathbb{C}^2)^2$) relation
\[ (e_{10}^2-4e_{20})(e_{01}^2-4e_{02}) = (e_{10}e_{01}-2e_{11}). \]

The mapping $\mathcal{A}$ is a birational equivalence, thus
$\mathcal{A}^*$ induces isomorphism of the field of Abelian
functions  $\mathcal{F}(\Jac V)$ on $\Jac V$ with the field of
rational functions $\mathcal{F}\left((V)^2\right)$ on $(V)^2$. We
have:
\[ \mathcal{A}^*(\wp_{11}(u)) = e_{10},\quad  \mathcal{A}^*(\wp_{13}(u)) =
-e_{20},\quad \mathcal{A}^*(\wp_{111}(u)) =
\frac{e_{10}e_{01}-2e_{11}}{e_{10}^2-4e_{20}}.  \]

In this way, the Theorem  \ref{T-1} completely describes the
isomorphism $\mathcal{A}^*$ and gives an opportunity to express
explicitly, for example, the even functions
$\wp_{klp}(u)\wp_{k'l'p'}(u)$ as polynomials of $\wp_{kl}(u)$. The
explicit formulae for those polynomials are given in the book
\cite{BEL}.

The hyperelliptic involution acts on $(V)^2$ according to the
formula
\[ \left[( s_1,\mu_1), (s_2,\mu_2) \right] \longrightarrow \left[( s_1,-\mu_1), (s_2,-\mu_2) \right]. \]
The Abel mapping is invariant with respect to this involution on
$(V)^2$ and the involution   $u \to -u$\, on \, $\Jac V$.

From the above formulae one can see that the images of the functions
$\wp_{kl}(u)$ and $\wp_{klp}(u)\wp_{k'l'p'}(u)$ are even, while the
images of the functions  $\wp_{klp}(u)$ are odd. Thus, the mapping
\[ \widehat{\mathcal{A}} \colon (V)^2/_{\pm} \longrightarrow K=(\Jac V)/_{\pm}, \]
is defined and it induces a homomorphism between rings of functions.

There is an addition law on $(V)^2$, with the Abel mapping
$\mathcal{A}$  as a homomorphism. Explicit form of this operation in
the coordinates  $\left[( s_1,\mu_1), (s_2,\mu_2) \right]$ has been
described in \cite{BL}.

Thus, on $(V)^2/_{\pm}$ there is corresponding two-valued addition,
such that   $\widehat{\mathcal{A}}$ is a homomorphism with respect
to the two-valued group structure on the Kummer variety $K$, defined
above.

\medskip

\

\section{Solutions of the system of equations of S. V. Kowalevski in genus two
$\wp$-functions}\label{sec:solkow}
\medskip

We follow chapter IV of the Golubev book \cite{Gol}.

Kowalevski introduced  variables $s_1,\; s_2$, which satisfy the
system of equations (see equations (17) and(18) from \cite{Gol}:
\begin{equation}\label{K-1}
\aligned \frac{ds_1}{\sqrt{\Phi(s_1)}} +
\frac{ds_2}{\sqrt{\Phi(s_2)}} &= 0,
\endaligned\end{equation}
\begin{equation}
\aligned
\frac{s_1ds_1}{\sqrt{\Phi(s_1)}} +
\frac{s_2ds_2}{\sqrt{\Phi(s_2)}} &= \frac{i}{2}\,dt,\label{K-2}
\endaligned
\end{equation}
where $\Phi(s)$ is a polynomial of fifth degree. The Abel mapping
\[ \mathcal{A} \colon (V)^2 \longrightarrow \Jac V \]
where $V = \left\{(s,\mu) \in \mathbb{C}^2\;:\; \mu^2 = \Phi(s)
\right\}$ è $\Phi(s) = s^5+\lambda_4s^3+\lambda_6s^2
+\lambda_8s+\lambda_{10}$, is defined by the system of equations
\begin{align}\label{A-3}
\frac{ds_1}{\mu_1} + \frac{ds_2}{\mu_2} &= du_3,\\
\frac{s_1ds_1}{\mu_1} + \frac{s_2ds_2}{\mu_2} &= du_1.\label{A-4}
\end{align}
We are going to describe solutions of the system of equations
(\ref{K-1}), (\ref{K-2}), following the book \cite{BEL}.

Let us consider the sigma-function $\sigma(u) = \sigma(u,\lambda)$,
associated with the curve $V$, and \linebreak corresponding Abelian
functions
$$\wp_{k,l}(u) = -\frac{\partial}{\partial u_k
\partial u_l}\ln \sigma(u),\; k=1,3.
$$
The general solution of the system of equations
(\ref{A-3}), (\ref{A-4}), which is the solution of the Jacobi
inversion problem for the Abel mapping, is represented by the point
$\left[ (s_1,\mu_1),(s_2,\mu_2) \right]$, where $(s_1,s_2)$ are the
solutions of the equation
\begin{equation}\label{A-5}
s^2-\wp_{11}(u)s-\wp_{13}(u) = 0,
\end{equation}
and for given $s_l$ we have
\begin{equation}\label{A-6}
2\mu_l = \wp_{111}(u)s_l + \wp_{113}(u), \; l=1,2.
\end{equation}
Thus, we get the solution  $\left[ \big( s_1(t),\mu_1(t) \big),
\big( s_2(t),\mu_2(t) \big) \right] \in (V)^2$ of the Kowalevski
system of equations (\ref{K-1}), (\ref{K-2}) in terms of genus 2
$\wp$-functions, where $\big( s_l(t),\mu_l(t) \big) = \big(
s_l(t,c_1,c_3),\mu_l(t,c_1,c_3) \big),\; l=1,2$, is the solution of
the Jacobi problem for  $u=(u_1,u_3)\in \Jac V$, where $u_1 = u_1(t)
= c_1+it,\; u_3 = c_3$  and  $c_1,\; c_3$ are constants.

The solution of the system (\ref{K-1}), (\ref{K-2}) based on the
sigma - functions, unlike the classical solutions in terms of the
theta-functions, is stable under the limit procedure when  $\lambda
\to 0$.

Without changing the notation, let us consider corresponding
functions, associated with the curve $V_0 = \left\{(s,\mu) \in
\mathbb{C}^2\;:\; \mu^2 = s^5 \right\}$. We have
\begin{align*}
\sigma(u) &= u_3-{1\over 3}\,u_1^3,\\
\sigma(u)^2\wp_{11}(u) &= 2u_1\left(u_3+{1\over 6}\,u_1^3\right),\\
\sigma(u)^2\wp_{13}(u) &= - u_1^2.
\end{align*}
Thus, for $\sigma(u) \neq 0$, the equation (\ref{A-5}) with $\lambda
\to 0$ is equivalent to
\begin{equation}\label{A-7}
\widehat s\,^2 - 2u_1\left(u_3+{1\over 6}\,u_1^3\right)\widehat s +
u_1^2 \sigma(u)^2 = 0,
\end{equation}
where $\widehat s = \sigma(u)^2s$. Further,
\begin{align*}
\sigma(u)^3\wp_{111}(u) &= \sigma(u)^3 \frac{\partial}{\partial
u_1}\,
\wp_{11}(u) =  2\left[u_3^2+{7\over 3}\,u_1^3u_3+{1\over 9}\,u_1^6\right],\\
\sigma(u)^3\wp_{113}(u) &= \sigma(u)^3 \frac{\partial}{\partial
u_1}\,\wp_{13}(u) = - 2u_1\left(u_3+{2\over 3}\,u_1^3\right).
\end{align*}
In this way, for  $\sigma(u) \neq 0$, the formula (\ref{A-6}) with
$\lambda \to 0$ is equivalent to
\begin{equation}\label{A-8}
\widehat \mu_l = \left[u_3^2+{7\over 3}\,u_1^3u_3+{1\over
9}\,u_1^6\right]\widehat s_l - u_1\left(u_3+{2\over
3}\,u_1^3\right)\sigma(u)^2,
\end{equation}
where $\widehat \mu_l = \sigma(u)^5\mu_l$. Observe, that
$\widehat\mu_l^2 = \widehat s\,^5$.

Thus, the general solution of the inversion problem $\left[(
s_1,\mu_1), ( s_2,\mu_2) \right] \in (V_0)^2$, where $V_0 =
\left\{(s,\mu) \in \mathbb{C}^2\;:\; \mu^2 = s^5 \right\}$, has the
form
\begin{equation}\label{A-9}
\sigma(u)^2 s_{1,2} = u_1\left[\left(u_3+{1\over 6}\,u_1^3\right)
\pm \sqrt{u_1^3\left(u_3-{1\over {12}}\,u_1^3\right)}\,\right]\,,
\end{equation}
and $\widehat \mu_l = \sigma(u)^5 \mu_l$ are given by the formula
(\ref{A-8}).

Let us consider important particular cases:
\begin{enumerate}
\item Let $u_1=0$, then $\sigma(u)=u_3$ and
\[ s_1=s_2=0,\quad \mu_1=\mu_2=0. \]
\item Let $u_3=0$, then $\sigma(u)=-{1 \over 3}\,u_1^3$ and
\[ \widehat s_{1,2} = {1 \over 6}\left(1\pm i\sqrt{3}\right)u_1^4, \quad \widehat \mu_{1,2}
= {1 \over {18}}\left(-1\pm \frac{i}{\sqrt{3}}\right)u_1^{10}.  \]
Check: the identity holds
\[ {1 \over {18^2}}\left(-1 + \frac{i}{\sqrt{3}}\right)^2 = {1 \over {6^5}}\left(1+ i\sqrt{3}\right)^5. \]
By taking  $u_1=c_1+\frac{i}{2}t$ and $u_3=0$, we get particular
solution of the Kowalevski system in the rational limit.
\item Let $u_3+{2 \over 3}\,u_1^3=0$, then $\sigma(u)=u_1^3$ è
\[ \widehat s_{1,2} = {1 \over 2}\left(-1\pm i\sqrt{3}\right)u_1^4, \quad \widehat
\mu_l = -\widehat s_lu_1^6. \] Check: the identity holds
$\left(-1\pm i\sqrt{3}\right)^3 = 8$.
\item Let $u_3 = {1 \over {12}}\,u_1^3$, then $\sigma(u)=-{1 \over 4}u_1^3$  è
\[ \widehat s_{1,2} = {1 \over {2^2}}\,u_1^4, \quad \widehat
\mu_{1,2} = {1 \over {2^5}}\,u_1^{10}. \]

\end{enumerate}
\medskip

\

\

\section{Geometric two-valued group laws and Kummer
varieties}\label{sec:kummergeom}

\

\subsection{A quadric in $\mathbb CP^5$ and a line complex in $\mathbb CP^3$}

\

Following classics, let us consider a three-dimensional projective
space $\mathbb CP^3= P(\mathbb C^4)$ and corresponding Grassmannian
$Gr(2,4)$ of all lines in $\mathbb CP^3$. By Pl\"ucker embedding,
the Grassmannian $Gr(2,4)$ can be realized as a quadric $G$ in
$\mathbb CP^5= P(\wedge^2 V)$, where $V=\mathbb C^4$:
$$
\aligned
Gr(2,4)&\hookrightarrow \mathbb CP^5\\
\ell=\ell<v_1, v_2>&\mapsto v_1 \wedge v_2, \quad v_1, v_2 \in V^4.
\endaligned
$$
The quadric $G$ parameterizes all decomposable elements $w=v_1\wedge
v_2$ of $P(\wedge^2 V)$ and the quadric is described by the Plucker
quadratic relation:
$$
G: w\wedge w =0.
$$
For a given element $x\in G$, denote by $\ell_x$ the line in
$\mathbb CP^3$ which maps to $x$ by the above embedding.

We consider so-called Schubert cycles:
$$
\aligned \sigma_1(\ell)&=\{x\in G\mid \ell_x \cap\ell \neq
\emptyset\}\\
\sigma_2(p)&=\{x\in G\mid \ell_x \ni p \}\\
\sigma_{1,1}(h)&=\{x\in G\mid \ell_x \subset h \}\\
\sigma_{2,1}(p,h)&=\{x\in G\mid p\in \ell_x \subset h \}\\
\endaligned
$$
with the intersection table
$$
\aligned
\sigma_1\cdot \sigma_1&=\sigma_2+\sigma_{1,1}\\
\sigma_1\cdot \sigma_2&=\sigma_1\cdot \sigma_{1,1}=\sigma_{2,1}\\
\sigma_2\cdot \sigma_2&=\sigma_{1,1}\cdot\sigma_{1,1}=\sigma_1\cdot
\sigma_{2,1}=1\\
\sigma_2\cdot \sigma_{1,1}&=0.
\endaligned
$$
One can easily see that every cycle $\sigma_1(\ell)$ is a hyperplane
section of the quadric $G$. If the line $\ell \in \mathbb CP^3$ is
determined by vectors $v_1, v_2\in V^4$ then the hyperplane of
intersection is of the form $H_{v_1\wedge v_2}=\{w\mid w\wedge v_1
\wedge v_2=0\}$.

Every cycle $\sigma_{2,1}(p,h)$ is a line in $\mathbb CP^5$. Every
line $L\subset G$ is of the form $L=\sigma_{2,1}(p,h)$.

A line $L\subset G$ is a pencil of lines in $\mathbb CP^3$. This is
a confocal pencil with the common point, {\it the focus} $p\in
\mathbb CP^3$. At the same time this pencil is coplanar with the
common plane $h\in \mathbb CP^{3*}$.

Every cycle of the form $\sigma_2(p)$ or of the form
$\sigma_{1,1}(h)$ is a two-plane in $\mathbb CP^5$. Conversely,
every two-plane in $G$ is of the form $\sigma_2(p)$ or of the form
$\sigma_{1,1}(h)$.

\

Let us recall some general properties of a quadric $Q$ in $\mathbb
CP^m$. The rank of the quadric is equal to the rank of any of its
symmetric $(m+1)\times (m+1)$ matrices. The quadric is smooth if its
rank is maximal, i.e. if rank of $Q$ is $m+1$. If the rank of $Q$ is
$r$ then it is a cone over a smooth quadric in $\mathbb CP^{r-1}$
with a vertex $\mathbb CP^{m-r}$. Quadrics of rank $m+1$ and $m$ are
called {\it general}.

\medskip

\begin{lemma}\label{lemma:quadrics}
Let $Q\in \mathbb CP^{2m+1}$ be a general quadric.
\begin{itemize}
\item[(a)] The dimension of a maximal linear subspace of $Q$ is equal
to $m$.
\item[(b)] The collection of maximal linear subspaces $C(Q)$ forms an
algebraic variety of dimension $m(m+1)/2$.
\item[(c)] If the rank of $Q$ is even, then $C(Q)$ has two components.
Otherwise, the component is unique. A component is unirational
variety.
\item[(d)] Let $L\subset Q$ be a linear subspace of dimension $m-1$
which does not contain the vertex of $Q$. For an irreducible family
of maximal subspaces $A\in C(Q)$, there is a unique maximal subspace
$M=M(L,A)$ of dimension $m$ which belongs to $A$ and contains $L$.
\item[(e)] Suppose that $Q$ is smooth and let $M_1, M_2$ be two
of its maximal linear subspaces. Then
$$
\dim M_1\cap M_2 \equiv m \,(\mod 2)
$$
is equivalent to the fact that $M_1, M_2$ belong to the same
irreducible component of $C(Q)$.
\end{itemize}
\end{lemma}

\medskip

Now, we specialize previous statement for the case of a smooth
quadric in $\mathbb CP^5$.

\medskip

\begin{proposition}\label{prop:quadric}
Let $G$ be a smooth quadric in $\mathbb CP^5$. Then:
\begin{itemize}
\item[(a)] There exists a four-dimensional vector space $V^4$ such
that $G$ is the Plucker embedding of Grassmannian of lines in $
P(V^4)$.
\item[(b)] Maximal linear subspaces of $G$ are two-dimensional and
they are of the form $\sigma_2(p)$ or of the form $\sigma_{1,1}(h)$,
where $p\in \mathbb CP^3= P(V^4)$ and $h\in \mathbb CP^{3*}$.
\item[(c)] Variety $C(Q)$ of all two-dimensional subspaces of $G$ is
three-dimensional and it has two irreducible components, $A$ and
$B$:
$$
\aligned A&=\{\sigma_2(p)\mid p\in \mathbb CP^3= P(V^4)\}\\
B&=\{\sigma_{1,1}(h)\mid h\in \mathbb CP^{3*}\}.
\endaligned
$$
\end{itemize}
\end{proposition}

\medskip

\subsection{Smooth intersection of two quadrics and Abelian
varieties}

\medskip

Now we consider the intersection $X$ of two quadrics $G$ and $F$ in
$\mathbb CP^5$. Such a set is classically called quadratic complex
of lines, if $G$ is understood as a Grassmannian of lines in some
$\mathbb CP^3$. Together with two quadrics $G$ and $F$, one may
consider the whole {\it pencil of quadrics}:
$$
F_{\lambda}:= F +\lambda G,
$$
and $X$ is the base set for the pencil, the common intersection of
quadrics from the pencil.

A pencil of quadrics is {\it generic} if associated pencil of
$6\times 6$ symmetric matrices contains six different singular
matrices. For a generic pencil $F_{\lambda}$ denote by $\lambda_1,
\dots, \lambda_6$ corresponding values of the pencil parameter
associated with singular matrices.

The condition that $X$ is smooth is equivalent to the condition that
the pencil $F_{\lambda}$ is generic. Smoothness of $X$ is also
equivalent to the fact that all quadrics $F_{\lambda}$ are general
and that exactly six of them, $F_{\lambda_1},\dots, F_{\lambda_6}$
are singular.
\medskip
\begin{proposition}\label{prop:X}
Suppose $X$ is smooth intersection of two quadrics $X=G\cap F$ in
$\mathbb CP^5$. Then:
\begin{itemize}
\item[(a)] Maximal linear subspaces of $X$ are one-dimensional.
\item[(b)] There is a maximal linear subspace through each point
of $X$.
\item[(c)] There are four maximal linear subspaces passing through a
generic point of $X$.
\end{itemize}
\end{proposition}

\medskip

Given a smooth intersection of quadrics $X=G\cap F$ in $\mathbb
CP^5$, following Narasimhan, Ramanan, Reid and Donagi, let us
consider the set of one-dimensional linear subspaces
$$
\mathcal A(X)=\{L\mid L\in X\cap Gr(2,6)\}.
$$
Suppose $G$ is realized as a Grassmannian of lines in some $\mathbb
CP^3$, and denote as before the two components of $C(G)$ of
two-dimensional linear subspaces of $G$ as
$$
\aligned A&=\{\sigma_2(p)\mid p\in \mathbb CP^3=P(V^4)\}\\
B&=\{\sigma_{1,1}(h)\mid h\in \mathbb P^{3*}\}.
\endaligned
$$

\medskip
\begin{lemma}\label{lemma:involutions}
Given $L\in \mathcal A(X)$. Then:
\begin{itemize}
\item[(a)] There is a unique two-dimensional linear subspace of $G$,
$\sigma_2(p)\in A$ such that $L\subset \sigma_2(p)$. There is a
unique one-dimensional linear subspace of $X$
$$
L_1\in \mathcal A(X),
$$
such that
$$
\sigma_2(p)\cap F=L\cup L_1.
$$
\item[(b)] There is a unique two-dimensional linear subspace of $G$,
$\sigma_{1,1}(h)\in A$ such that $L\subset \sigma_{1,1}(h)$. There
is a unique one-dimensional linear subspace of $X$
$$
L_2\in \mathcal A(X),
$$
such that
$$
\sigma_{1,1}(h)\cap F=L\cup L_2.
$$
\end{itemize}
\end{lemma}
\medskip
The last Lemma introduces two involutions
$$
\aligned i_1:&\mathcal A(X)\rightarrow \mathcal A(X)\\
i_1:&L\mapsto L_1\\
i_2:&\mathcal A(X)\rightarrow \mathcal A(X)\\
i_2:&L\mapsto L_2.
\endaligned
$$
As pencils of lines in $\mathbb CP^3$ the subspaces $L$ and $L_1$
are confocal, they have the same focus $p$. In the same manner, the
subspaces $L$ and $L_2$ are coplanar, they have the same plane $h$.

Moreover, there are two mappings
$$
\aligned k_1:&\mathcal A(X)\rightarrow \mathbb CP^3\\
k_1:&L\mapsto p\\
k_2:&\mathcal A(X)\rightarrow \mathbb CP^{3*}\\
k_2:&L\mapsto h.
\endaligned
$$
The mapping $k_1$ maps a pencil $L$ to its focus in $\mathbb CP^3$
while $k_2$ maps a pencil to its plane in $\mathbb CP^{3*}$.

Denote by $K\subset \mathbb CP^3$ the image of $\mathcal A(X)$ by
$k_1$. We see that $k_1$ is a double covering of $\mathcal A(X)$
over $K$ and that the involution $i_1$ interchanges the leaves of
the covering. We are going to call $K$ {\it the Kummer variety} of
$\mathcal A(X)$. It is associated to the choice of a quadric $G$
from the pencil and to the choice of a connected component $A$ of
$C(G)$.

Similarly, denote by $K^*\subset \mathbb CP^{3*}$ the image of
$\mathcal A(X)$ by $k_2$. The mapping $k_2$ is a double covering of
$\mathcal A(X)$ over $K$ and  the involution $i_2$ interchanges the
leaves of the covering. We are going to call $K^*$ {\it the dual
Kummer variety} of $\mathcal A(X)$. It is associated to the choice
of a quadric $G$ from the pencil and to the choice of a connected
component $B$ of $C(G)$.

It can be shown that $K^*$ is dual to $K$, which means that every
plane from $K^*$ is tangent to $K$. From the previous considerations
it can also be seen that the degree of $K$ is equal to 4.

For a general point $p\in \mathbb CP^3$ the plane
$\sigma_2(p)\subset G$ intersects $F$ along a smooth conic. The
Kummer variety can be described as a set of the points in $\mathbb
CP^3$ for which this intersection is a conic which is not smooth.
For general point of $K$ this intersection is a degenerate conic
which is a union of two different lines $L$ and $i_1(L)$. But, there
is a subset $R\in K$ of sixteen points, for which this intersection
is a degenerate conic of double line. These sixteen points from $R$
correspond to the fixed points of the involution $i_1$ and to the
ramification points of the double covering. We also denote
$R^*\subset K^*$ the set of points that correspond to the fixed
points of the involution $i_2$.

\medskip

\subsection{Pencils of quadrics, hyperelliptic curves and geometric group laws}

\medskip

With a generic pencil of quadrics $F_{\lambda}$ in $\mathbb CP^5$
with six singular quadrics $F_{\lambda_1},\dots, F_{\lambda_6}$ one
may associate a genus two curve
$$
\Gamma: y^2=\prod_{i=1}^6(x-\lambda_i).
$$
As it was shown by Reid and Donagi, this correspondence between a
generic pencil of quadrics and a hyperelliptic curve is not just
formal. After Donagi, denote by $E$ the family of all connected
components of $C(F_{\lambda})$ of all quadrics from the pencil
together with a projection
$$
p:E\rightarrow \mathbb CP^1
$$
which maps a given irreducible family of two-subspaces of a quadric
$F_{\mu}$ to the value $\mu$ of the pencil parameter. The projection
$p$ is obviously double covering. It is  ramified over the points
$\lambda_i, i=1, \dots, 6$ since the singular quadrics are the only
one with unique component of maximal linear subspaces. Thus, Donagi
showed the isomorphism between $E$ and $\Gamma$.

But, there is yet another natural realization of the hyperelliptic
curve $\Gamma$ in the context of the pencil $F_{\lambda}$, as it was
demonstrated by Reid.

For $L\in \mathcal A(X)$ denote by $\mathcal A_L(X)$ the closure of
the set $\{L'\in \mathcal A(X)\mid L\cap L' \ne \emptyset\}$. There
is a natural projection
$$q:\mathcal A_L(X)\setminus \{L\} \rightarrow \mathbb CP^1$$
which maps $L'$ to the parameter $\mu$ of a quadric $F_{\mu}$ if the
space $<L, L'>$ belongs to $C(F_{\mu})$. The mapping $q$  is double
covering ramified over the six points $\lambda_i, i=1, \dots, 6$ and
$\mathcal A_L(X)$ is isomorphic to the hyperelliptic curve $\Gamma$.
The natural involution on $\mathcal A_L(X)$ which interchanges the
folds of $q$ will be denoted by $\tau_L$.

Moreover, as it was shown by Reid and Donagi, $\mathcal A(X)$ is an
Abelian variety, isomorphic to the Jacobian of the curve $\Gamma$.

\medskip

We will refer to the curves $\mathcal A_L(X)$ and $E$ as {\it
Donagi-Reid-Kn\"orrer curves} (DRK) associated with intersection of
quadrics $X$ of a generic pencil.

\medskip

 It can easily be shown that for any hyperelliptic curve
$\Gamma$, there exists a pencil of quadrics with the base set $X$
such that $\mathcal A(X)$ is isomorphic to the Jacobian of $\Gamma$.

The addition laws on the Abelian varieties $\mathcal A(X)$ are such
that
$$L_1+L_2=M_1+M_2$$
if there exists $\mu$ such that
$$ <L_1, L_2>, <M_1,M_2>\in C(F_{\mu})$$
and if the two two-dimensional spaces $ <L_1, L_2>, <M_1,M_2>$
belong to the same component of $C(F_{\mu})$.

\medskip

\begin{proposition}
Every point $e\in E$ determines its Kummer variety involution $i_e$
and its Kummer variety $K_e$. The hyperelliptic involution on $E$,
$\tau$ interchanges a Kummer variety and its dual:
$$
K_{\tau(e)}=K_e^*.
$$
\end{proposition}

\medskip

We will fix a point $e_0\in E$ and a line $L_0\in \mathcal A(X)$ as
the origin of a group structure on $\mathcal A(X)$ such that
$$
L_0+\hat L_0=0,
$$
whenever
$$<L_0,\hat L_0>\in e_0.$$

As above, denote by $\tau_{L_0}$ the natural involution on $\mathcal
A_{L_0}(X)$. Then we have
\medskip
\begin{lemma}\label{lemma:involutions} The involutions $\tau_{L_0}$
and $i_{e_0}$ are related according to the formula
$$
i_{e_0}|_{\mathcal A_{L_0}(X)}=\tau_{L_0}.
$$
\end{lemma}
\medskip

Now, we are ready to define a two-valued group structure on the
Kummer variety $K_{e_0}$. We will define a mapping
$$
\star: K_{e_0}\times K_{e_0}\longrightarrow (K_{e_0})^2
$$
by the following procedure.

Take a pair
$$(p_1, p_2)\in K_{e_0}\times K_{e_0}.$$
Denote by
$$ [L_1]=\{L_1,\hat L_1\}, \quad L_1, \hat L_1\in \mathcal A(X)$$
the class of lines in $\mathcal A(X)$ which represents the two
confocal pencils of lines in $\mathbb CP^3$ with the focal point
$p_1$. Similarly, denote by
$$ [L_2]=\{L_2,\hat L_2\}, \quad L_2, \hat L_2\in \mathcal A(X)$$
the class of lines in $\mathcal A(X)$ which represents the two
confocal pencils of lines in $\mathbb CP^3$ with the focal point
$p_2$.

Assume that $L_1, L_2$ don't intersect $L_0$. Denote by $N_1, N_2$
the two lines in $\mathcal A_{L_0}(X)$ of intersection of the space
$<L_0,L_1>$ with $X$ and by $e_1, e_2\in E$ denote the classes
determined by
$$
<L_0,N_1>\in e_1, \quad <L_0, N_2>\in e_2.
$$
Denote also by
$$\mu_1=p(e_1), \quad \mu_2=p(e_2),$$
and by
$$N_1',N_1'', N_2',N_2''\in \mathcal A_{L_2}(X)$$
the lines which intersect $L_2$ and which are uniquely defined by
the conditions
$$
\aligned &<N_1', L_2>\in e_1, \quad <N_1'', L_2>\in \tau(e_1)\\
&<N_2', L_2>\in e_1, \quad <N_2'', L_2>\in \tau(e_2).
\endaligned
$$
In other words $N_1', N_1''$ belong to the two two-dimensional
spaces of the different classes of $C(F_{\mu_1})$ which contain
$L_2$; $N_2', N_2''$ belong to the two two-dimensional spaces of the
different classes of $C(F_{\mu_2})$ which contain $L_2$.

Let $M_1$ be the fourth intersection line in the intersection of $X$
with the space generated with $L_2, N_1', N_2'$ and let $M_2$ be the
fourth intersection line in the intersection of $X$ with the space
generated with $L_2, N_1'', N_2''$. The line $M_1$ represents a
pencil of lines in $\mathbb CP^3$ with the focal point $w_1$ and
$M_2$ represents a pencil of lines in $\mathbb CP^3$ with the focal
point $w_2$.

If we repeat the above procedure with $\hat L_1, \hat L_2$ instead
of $L_1, L_2$ we come to the lines $\hat M_1$ and $\hat M_2$ which
are confocal with $M_1$ and $M_2$.

One can easily adjust previous construction to the case where $L_1,
L_2$ intersect $L_0$: denote by $M_1$ the fourth line of the
intersection of $X$ with the space generated with $\hat L_1, \hat
L_2, L_0$; denote by $M_2$ the fourth line of the intersection of
$X$ with the space generated with $L_1, \hat L_2, L_0$.

Thus we get
\medskip

\begin{theorem}\label{th:geom2group}
The mapping defined by the formulae
$$
\aligned
 \star: &K_{e_0}\times K_{e_0}\longrightarrow (K_{e_0})^2\\
&p_1\star p_2 = (w_1,\,w_2)
\endaligned
$$
defines a structure of two-valued group on the Kummer variety
$K_{e_0}$.
\end{theorem}
\medskip

\begin{definition} The mapping $\star$ defined by previous
construction defines {\it the geometric two-valued group law} on the
Kummer variety $K_{e_0}$. \end{definition}

\section{Integrable billiards and two-valued  laws}\label{sec:int}

\medskip

\subsection{Pencils of quadrics and billiards, an overview}

We begin this Section by repeating basic definitions related to
billiard systems of confocal quadrics from \cite {DR}, \cite {DR2}.

\smallskip

Let $\mathcal Q_1$ and $\mathcal Q_2$ be two quadrics. Denote by $u$
the tangent plane to $\mathcal Q_1$ at point $x$ and by $z$ the pole
of $u$ with respect to $\mathcal Q_2$. Suppose lines $\ell_1$ and
$\ell_2$ intersect at $x$, and the plane containing these two lines
meet $u$ along $\ell$.

\begin{definition}{\rm If lines $\ell_1, \ell_2, xz, \ell$ are
coplanar and harmonically conjugated, we say that rays $\ell_1$ and
$\ell_2$ {\it obey the reflection law} at the point $x$ of the
quadric $\mathcal Q_1$ with respect to the confocal system which
contains $\mathcal Q_1$ and $\mathcal Q_2$.} \end{definition}

If we introduce a coordinate system in which quadrics $\mathcal Q_1$
and $\mathcal Q_2$ are confocal in the usual sense, reflection
defined in this way is same as the standard one.

\begin{theorem}[One Reflection Theorem]
Suppose rays $\ell_1$ and $\ell_2$ obey the reflection law at $x$ of
$\mathcal Q_1$ with respect to the confocal system determined by
quadrics $\mathcal Q_1$ and $\mathcal Q_2$. Let $\ell_1$ intersects
$\mathcal Q_2$ at $y_1'$ and $y_1$, $u$ is a tangent plane to
$\mathcal Q_1$ at $x$, and $z$ its pole with respect to $\mathcal
Q_2$. Then lines $y_1'z$ and $y_1z$ respectively contain
intersecting points $y_2'$ and $y_2$ of ray $\ell_2$ with $\mathcal
Q_2$. Converse is also true.
\end{theorem}

\begin{corollary} Let rays $\ell_1$ and $\ell_2$ obey the
reflection law of $\mathcal Q_1$ with respect to the confocal system
determined by quadrics $\mathcal Q_1$ and $\mathcal Q_2$. Then
$\ell_1$ is tangent to $\mathcal Q_2$ if and only if is tangent
$\ell_2$ to $\mathcal Q_2$; $\ell_1$ intersects $\mathcal Q_2$ at
two points if and only if $\ell_2$ intersects $\mathcal Q_2$ at two
points.
\end{corollary}

Next assertion is crucial for applications to the billiard dynamics.

\begin{theorem} [Double Reflection Theorem] \label{th:DRT}
Suppose that $\mathcal Q_1$, $\mathcal Q_2$ are given quadrics and
$\ell_1$ line intersecting $\mathcal Q_1$ at the point $x_1$ and
$\mathcal Q_2$ at $y_1$. Let $u_1$, $v_1$ be tangent planes to
$\mathcal Q_1$, $\mathcal Q_2$ at points $x_1$, $y_1$ respectively,
and $z_1$, $w_1$ their poles with respect to $\mathcal Q_2$ and
$\mathcal Q_1$. Denote by $x_2$ second intersecting point of the
line $w_1x_1$ with $\mathcal Q_1$, by $y_2$ intersection of $y_1z_1$
with $\mathcal Q_2$ and by $\ell_2$, $\ell_1'$, $\ell_2'$ lines
$x_1y_2$, $y_1x_2$, $x_2y_2$. Then pairs $\ell_1,\ell_2$;
$\ell_1,\ell_1'$; $\ell_2,\ell_2'$; $\ell_1',\ell_2'$ obey the
reflection law at points $x_1$ (of $\mathcal Q_1$), $y_1$ (of
$\mathcal Q_2$), $y_2$ (of $\mathcal Q_2$), $x_2$ (of $\mathcal
Q_1$) respectively.
\end{theorem}

\begin{corollary}
If the line $\ell_1$ is tangent to a quadric $\mathcal Q$ confocal
with $\mathcal Q_1$ and $\mathcal Q_2$, then rays $\ell_2$,
$\ell_1'$, $\ell_2'$ also touch $\mathcal Q$.
\end{corollary}

\medskip
For the conclusion, we recall the notion of generalized Cayley's
curve from \cite {DR}, \cite {DR2}.

\begin{definition}{\rm
The {\it generalized Cayley curve} $\mathcal C_{\ell}$ is the
variety of hyperplanes tangent to quadrics of a given confocal
family in $\mathbb {C}P^d$ at the points of a given line $\ell$. }
\end{definition}

This curve is naturally embedded in the dual space $\mathbb
CP^{d\,*}$.

\medskip
\begin{proposition} The generalized Cayley curve in
$\mathbb {C}P^d$, for $d\ge 3$ is a hyperelliptic curve of genus
$g=d-1$. Its natural realization in $\mathbb {C}P^{d\,*}$ is of
degree $2d-1$.
\end{proposition}

The natural involution $\tau_{\ell}$ on the generalized Cayley's
curve $\mathcal C_{\ell}$ maps to each other the tangent planes at
the points of intersection of $\ell$ with any quadric of the
confocal family.

\smallskip

It was observed in \cite{DR} that this curve is isomorphic to the
Veselov-Moser isospectral curve.

\medskip
Now we are going to mention a connection obtained in \cite {DR2},
between generalized Cayley's curve defined above and the curves (see
the previous Section) studied by Kn\"orrer, Donagi, Reid. This
connetion  traces out the relationship between billiard
constructions and the algebraic structure of the corresponding
Abelian varieties.

\smallskip

The famous Chasles theorem (see \cite {Ar}) states that any line in
the space $\mathbf R^d$ is tangent to exactly $d-1$ quadrics from a
given confocal family:
\begin{equation}\label{eq:konfokalna.familija}
\mathcal Q_{\lambda}:\ Q_{\lambda}(x) = 1.
\end{equation}
where we denote:
$$
Q_{\lambda}(x)=\frac {x_1^2}{a_1-\lambda}+\dots + \frac
{x_d^2}{a_d-\lambda}.
$$
We assume that the family (\ref{eq:konfokalna.familija}) is generic,
i.e.\ the constants $a_1,\dots,a_d$ are all distinct.

\smallskip

 Suppose a line $\ell$ is tangent to quadrics
 $\mathcal Q_{\alpha_1},\dots,\mathcal Q_{\alpha_{d-1}}$
from the given confocal family. Denote by $\mathcal A_{\ell}$ the
family of all lines which are tangent to the same $d-1$ quadrics.
Note that according to the corollary of the One Reflection Theorem,
the set $\mathcal A_{\ell}$ is invariant to the billiard reflection
on any of the confocal quadrics.

\smallskip

For what follows, the next simple observation is important.

\begin{lemma}\label{lema:veza}
Let the lines $\ell$ and $\ell'$ obey the reflection law at the
point $z$ of a quadric $\mathcal Q$ and suppose they are tangent to
a confocal quadric $\mathcal Q_1$ at the points $z_1$ and $z_2$.
Then the intersection of the tangent spaces $T_{z_1}\mathcal Q_1\cap
T_{z_2}\mathcal Q_1$ is contained in the tangent space
$T_{z}\mathcal Q$. \end{lemma}

\medskip
Following \cite{Kn}, together with $d-1$ affine confocal quadrics
$\mathcal Q_{\alpha_1},\dots,\mathcal Q_{\alpha_{d-1}}$, one can
consider their projective closures $\mathcal
Q^p_{\alpha_1},\dots,\mathcal Q^p_{\alpha_{d-1}}$ and the
intersection $X$ of two quadrics in $\mathbb {C}P^{2d-1}$:
\begin{equation}\label{eq:kvadrika1}
x_1^2+\dots+x_d^2-y_1^2-\dots - y_{d-1}^2=0,
\end{equation}
\begin{equation}\label{eq:kvadrika2}
 a_1x_1^2+\dots +a_dx_d^2 - \alpha_1y_1^2 - \dots -\alpha_{d-1}y_{d-1}^2=x_0^2.
\end{equation}

Denote by $\mathcal A(X)$ the set of all $(d-2)$-dimensional linear
subspaces of $X$. For a given $L\in \mathcal A$, denote by $\mathcal
A_L(X)$ the closure in $\mathcal A(X)$ of the set
 $\{\, L'\in F \mid \dim L\cap L' = d-3\, \}$.
It was shown in \cite{Re} that $F_L$ is a nonsingular hyperelliptic
curve of genus $d-1$.

\smallskip

The projection
$$
\pi'\ :\ \mathbb {C}P^{2d-1}\setminus\{(x,y)|x=0\}\to\mathbb {C}P^d,
\quad \pi'(x,y)=x,
$$
maps $L\in \mathcal A(X)$ to a subspace $\pi'(L)\subset\mathbb
{CP}^d$ of the codimension $2$. $\pi'(L)$ is tangent to the quadrics
$\mathcal Q^{p*}_{\alpha_1},\dots,\mathcal Q^{p*}_{\alpha_{d-1}}$
that are dual to $\mathcal Q^p_{\alpha_1},\dots,\mathcal
Q^p_{\alpha_{d-1}}$.

\smallskip

Thus, the space dual to $\pi'(L)$, denoted by $\pi^*(L)$, is a line
tangent to the quadrics $\mathcal Q^p_{\alpha_1},\dots,\mathcal
Q^p_{\alpha_{d-1}}$.

\smallskip

We can reinterpret the generalized Cayley's curve $\mathcal
C_{\ell}$, which is a family of tangent hyperplanes, as a set of
lines from $\mathcal A_{\ell}$ which intersect $\ell$. Namely, for
almost every tangent hyperplane there is a unique line $\ell'$,
obtained from $\ell$ by the billiard reflection. Having this
identification in mind, it is easy to prove the following

\begin{corollary} There is a birational morphism between the generalized
Cayley's curve $\mathcal C_{\ell}$ and Reid-Donagi-Kn\"orrer's curve
$F_L$, with $L=\pi^{*-1}(\ell)$, defined by
$$
j: \ell'\mapsto L',\quad L'=\pi^{*-1}(\ell'),
$$
where $\ell'$ is a line obtained from $\ell$ by the billiard
reflection on a confocal quadric.
\end{corollary}

\medskip
Lemma \ref{lema:veza}, giving a link between the dynamics of
ellipsoidal billiards and algebraic structure of certain Abelian
varieties, provides a two way interaction: to apply algebraic
methods in the study of the billiard motion, but also vice versa, to
use billiard constructions in order to get more effective, more
constructive and more observable understanding of the algebraic
structure.

\medskip

\subsection{Two-valued billiard structure}

\medskip

Thus, we are going to apply this relationship to construct a
billiard analogue of the geometric two-valued group structure
defined on the Kummer variety $K_{e_0}$ from the previous Section,
see Theorem \ref{th:geom2group} below.

\smallskip

In order to fit together notations from the previous Section and
from the last Subsection, we set $d=3$ and
$$
\{\lambda_1,\dots,\lambda_6\}=\{a_1,a_2,a_3,\alpha_1,\alpha_2,\infty\}.
$$
We consider the confocal pencil of quadrics $Q_{\lambda}$ in
$\mathbb CP^3$ and the set $\mathcal A_{\ell}$ of lines  in $\mathbb
CP^3$ which are tangent to both quadrics $Q_{\alpha_1},
Q_{\alpha_2}$.

We select a quadric $Q_{p(e_0)}$ together with a family of
reflections  which corresponds to $e_0$ and a line $\ell_0\in
\mathcal A_{\ell}$.

We introduce an involution $I_{e_0}$ on $\mathcal A_{\ell}$ induced
by $Q_{p(e_0)}$ with $e_0$,   which maps a line $m\in \mathcal
A_{\ell}$ to a line $m'\in \mathcal A_{\ell}$ uniquely defined by
the condition that $m$ and $m'$ are obtained by billiard reflection
of each other from $Q_{p(e_0)}$ from the given family $e_0$.

We suppose that the restriction of this involution $I_{e_0}$
restricted on $C_{\ell_0}$ coincides with the natural involution on
the generalized Cayley curve $C_{\ell_0}$.
\medskip

Denote by
$$ [\ell_1]=\{\ell_1,\hat \ell_1\}, \quad \ell_1, \hat \ell_1\in \mathcal A_{\ell}$$
and by
$$ [\ell_2]=\{\ell_2,\hat \ell_2\}, \quad \ell_2, \hat \ell_2\in \mathcal A_{\ell}$$
two classes of lines in $\mathcal A_{\ell}/I_{e_0}$.

Assume that $\ell_1, \ell_2$ don't intersect $\ell_0$. Denote by
$n_1, n_2$ the two lines in $C_{\ell_0}(X)$ which form a double
reflection configuration with $\ell_0, \ell_1$. Denote by
$Q_{\mu_1}$ the quadric of billiard reflection of
 $\ell_0,n_1$ and of billiard reflection of $\ell_1, n_2$. Similarly, denote by $Q_{\mu_2}$ the quadric of billiard
 reflections of pairs $\ell_0,n_2$ and $\ell_1,n_1$.

The quadrics $Q_{\mu_1}, Q_{\mu_2}$ intersect line $\ell_2$. Denote
 by
$$n_1',n_1'', n_2',n_2''\in  C_{\ell_2}$$
the lines of billiard  reflection of the line $\ell_2$ from the
quadrics $Q_{\mu_1}, Q_{\mu_2}$.

The lines  $n_1', n_1''$ are obtained by reflection at $Q_{\mu_1}$
and reflection of $n_1'$ belongs to the same family as the
reflection of $\ell_0, n_1$; the lines $n_2', n_2''$ are obtained by
reflection from $Q_{\mu_2}$ and reflection of $n_2'$ belongs to the
same family as the reflection $\ell_0,n_2$.

Let $m_1$ be the fourth  line of the double reflection configuration
determined with $\ell_2, n_1', n_2'$ and let $m_2$ be the fourth
line of the double reflection configuration generated by $\ell_2,
n_1'', n_2''$.

If we repeat the above procedure with $\hat \ell_1, \hat \ell_2$
instead of $\ell_1, \ell_2$ we come to the lines $\hat m_1$ and
$\hat m_2$ which are obtained from $m_1$ and $m_2$ by the involution
$I_{e_0}$.

One can easily adjust previous construction to the case where
$\ell_1, \ell_2$ intersect $\ell_0$: denote by $m_1$ the fourth line
of the double reflection configuration determined with $\hat \ell_1,
\hat \ell_2, \ell_0$; denote by $m_2$ the fourth line of the double
reflection configuration determined with $\ell_1, \hat \ell_2$ and $
\ell_0$.

Thus we get
\medskip

\begin{theorem}\label{th:geom2group}
The mapping defined by the formulae
$$
\aligned
 \star_b: &\mathcal A_{\ell}/I_{e_0}\times \mathcal A_{\ell}/I_{e_0}\longrightarrow (\mathcal A_{\ell}/I_{e_0})^2\\
&[\ell_1]\star_b [\ell_2] = ([m_1],\,[m_2])
\endaligned
$$
defines a structure of two-valued group on the  variety $\mathcal
A_{\ell}/I_{e_0}$.
\end{theorem}
\medskip

\begin{definition} The mapping $\star_b$ defined by previous
construction defines {\it the billiard two-valued group law} on the
 variety $\mathcal A_{\ell}/I_{e_0}$. \end{definition}

\medskip

\section{Moduli of semi-stable bundles and two-valued group
structure on Kummer varieties}\label{sec:moduli}

\medskip

As we have mentioned before, historically first examples of
$2$-valued groups appeared in topological context in the study of
the characteristic classes of vector bundles, see \cite {BN}. There
 one-dimensional symplectic bundles over $\mathbb {HP}^n$, and
together with the canonical projection $\mathbb
{C}P^{2n+1}\rightarrow \mathbb {H}P^n$, the associated
two-dimensional complex vector bundles over $\mathbb {C}P^{2n+1}$
with $c_1=0$, were considered.

For a pair of such two-dimensional bundles,
$$
\xi_1, \, \xi_2, \quad c_1(\xi_i)=0, \, i=1,2
$$
its tensor product
$$
\xi_1\otimes_C\xi_2
$$
is a four-dimensional bundle, with the first two Potryagin classes
$p_1(\xi_1\otimes_C\xi_2)$ and $p_2(\xi_1\otimes_C\xi_2)$. In
\cite{BN}, to the initial pair of bundles, a pair of {\it virtual}
two-bundles, were associated according to the formula
$$
Z^2-p_1(\xi_1\otimes_C\xi_2)Z+p_2(\xi_1\otimes_C\xi_2)=0.
$$
The solutions $Z_{1,2}$ of the last quadratic equation play a role
of the first Pontryagin classes of two virtual two-bundles. The last
quadratic equation defines a two-valued group structure in
$x=p_1(\xi_1)$ and $y=p_1(\xi_2)$, since
$p_1(\xi_1\otimes_C\xi_2)=\Theta_1(x, y)$ and
$p_2(\xi_1\otimes_C\xi_2)=\Theta_2(x, y)$, where $\Theta_i$ are
certain series.

In this Section, we are going to construct the two-valued group
structure on the Kummer variety, in the context of two-dimensional
varieties, semi-stable in the sense of Mumford and Seshadri. The
significance of the present situation, lies in the fact that
obtained resulting two-bundles {\it are not virtual} - they are
realized as a pair of two-bundles, semi-stable but not stable.

To get this, yet another interpretation of Kummer varieties, we are
going back to \cite{NR}  Let $X$ be a genus two curve. Following
Mumford and Seshadri, the notions of stable and semi-stable vector
bundles of rank $n$ and degree $d$ have been introduced
respectively. For a holomorphic nonzero vector bundle $W$ on $X$ one
introduces a rational-valued  function $\mu(W)=\deg W/ \rank W$. A
vector bundle $W$ is stable if for every proper subbundle $V$ the
condition
$$
\mu(V)<\mu(W)
$$
is satisfied. Similarly, a bundle is semi-stable if in the last
inequality the sign $<$ is replaced by $\le$. Any semi-stable bundle
$W$ has a strictly decreasing filtration
$$
W=W_0\supset W_1\supset\dots \supset W_n=\{0\}
$$
such that $W_{i-1}/W_{i}$ are stable and
$\mu(W_{i-1}/W_{i})=\mu(W)$. Denote by $Gr W=\oplus W_{i-1}/W_{i}$.
As Seshadri defined, two semi-stable bundles $W_1$ and $W_2$ are
S-equivalent if $Gr W_1\approx Gr W_2$; a normal, projective
$(n^2+1)$-dimensional variety of S-equivalence classes of
semi-stable bundles of degree d and rank n denote $U(n,d)$.

Following \cite{NR}, for $U(2,0)$, denote by $S$ its
three-dimensional sub-variety of bundles with trivial determinant.
The non-stable bundles in $S$ are of the form
$$
j\oplus j^{-1},
$$
with $j$ a line bundle of degree $0$. The Kummer surface $K$
associated to the Jacobian of $X$ is isomorphic to the set of all
non-stable bundles in $S$.

Now, we define a structure of two-valued group on $K$. It is an
important development of the Example \ref{ex:bundles} from Section
\ref{sec:def}. Denote $a, b \in K$, where
$$
a=j\oplus j^{-1}, \quad b=l\oplus l^{-1},
$$
where $j, l$ are line bundles on $X$ of degree $0$. Then:
\begin{equation}\label{eq:semistable}
a\star_s b:= (j\otimes l \oplus j^{-1}\otimes l^{-1}, j\otimes
l^{-1}\oplus j^{-1}\otimes l).
\end{equation}
\medskip

\begin{proposition} \label{prop:semistable}
The operation
$$
\star_s:K\times K \longrightarrow (K)^2
$$
determined by the relation \ref{eq:semistable} defines a two-valued
group structure on $K$.
\end{proposition}
\medskip

\subsection*{Acknowledgements}

The research of one of the authors (V. D.) was partially supported
by the Serbian Ministry of Science and Technological Development,
Project {\it Geometry and Topology of Manifolds, Classical Mechanics
and Integrable Dynamical Systems} and by the Mathematical Physics
Group of the University of Lisbon, Project \emph{ Probabilistic
approach to finite and infinite dimensional dynamical systems,
PTDC/MAT/104173/2008}.

\newpage\thispagestyle{empty}
\vspace*{20mm}

\end{document}